\documentclass{article}
\usepackage[english]{babel}

\usepackage[letterpaper,top=2cm,bottom=2cm,left=2cm,right=2cm,marginparwidth=1.75cm]{geometry}

\usepackage[mathscr]{euscript}
\usepackage[dvipsnames,table]{xcolor}
\usepackage[colorlinks=true, allcolors=OliveGreen]{hyperref}

\usepackage{calrsfs}
\usepackage{amssymb}
\usepackage{amsmath}
\usepackage{amsthm}
\usepackage{amsfonts}
\usepackage{subcaption}
\usepackage{multirow}
\usepackage{graphicx}
\usepackage{amsfonts}
\usepackage{mathtools}
\usepackage{mathrsfs}
\usepackage{placeins}
\usepackage{multirow}
\usepackage{verbatim,hyperref}
\usepackage{xcolor}
\usepackage{stmaryrd}
\usepackage{array}
\usepackage{caption}
\newcolumntype{M}[1]{>{\centering\arraybackslash}m{#1}}

\DeclareMathOperator{\diam}{diam}
\usepackage{algorithm}
\usepackage{algpseudocode}
\usepackage{dirtytalk}
\usepackage{todonotes}
\usepackage{authblk}
\usepackage{svg}
\usepackage{subcaption}
\usepackage{listings}
\usepackage{filecontents}
\usepackage{tikz}
\usepackage{pgfplots}

\newcolumntype{C}{>{\arraybackslash}p{7cm}}
\newcolumntype{Q}{>{\centering\arraybackslash}p{3.8cm}}
\newcolumntype{q}{>{\centering\arraybackslash}p{2.8cm}}
\newcommand{\jump}[1]{[\![#1]\!]}

\graphicspath{{Images/}}

\definecolor{codegreen}{rgb}{0,0.6,0}
\definecolor{codegray}{rgb}{0.5,0.5,0.5}
\definecolor{codepurple}{rgb}{0.58,0,0.82}
\definecolor{backcolour}{rgb}{0.95,0.95,0.92}

\title{Polytopal mesh agglomeration via geometrical deep learning for three-dimensional heterogeneous domains\footnote{\textbf{Funding}: This research has been funded by the European Union (ERC, NEMESIS, project number 101115663). Views and opinions expressed are, however, those of the author(s) only and do not necessarily reflect those of the European Union or the European Research Council Executive Agency. The present research is part of the activities of “Dipartimento di Eccellenza 2023-2027”. PFA and MC are members of INdAM-GNCS.}} 

\author[1]{Paola F. Antonietti\footnote{paola.antonietti@polimi.it}}

\author[1]{Mattia Corti\footnote{mattia.corti@polimi.it}}

\affil[1]{MOX-Dipartimento di Matematica, Politecnico di Milano, Piazza Leonardo da Vinci 32, Milan, 20133, Italy}

\author[ ]{Gabriele Martinelli\footnote{gabriele2.martinelli@mail.polimi.it}}

\begin{document}
\maketitle

\begin{abstract}
Agglomeration techniques can be successfully employed to reduce the computational costs of numerical simulations and stand at the basis of multilevel algebraic solvers. To automatically perform mesh agglomeration, we propose a novel Geometrical Deep Learning-based algorithm that can exploit the geometrical and physical information of the underlying computational domain to construct the agglomerated grid and -simultaneously-guarantee the agglomerated grid's quality. In particular, we propose a bisection model based on Graph Neural Networks (GNNs) to partition a suitable connectivity graph of computational three-dimensional meshes. The new approach has a high online inference speed. It can simultaneously process the graph structure of the mesh, the geometrical information of the mesh (e.g., elements' volumes, centers' coordinates), and the physical information of the domain (e.g., physical parameters). Taking advantage of this new approach, our algorithm can agglomerate meshes of a domain composed of heterogeneous media, automatically respecting the underlying heterogeneities. The proposed GNN approach is compared with the k-means algorithm and METIS, which are \textcolor{black}{widely employed}  approaches for graph partitioning and are meant to process only the connectivity information on the mesh. We demonstrate that the performance of our algorithms outperforms \textcolor{black}{the k-means and METIS algorithms} in terms of quality metrics and runtimes. Moreover, we demonstrate that our algorithm also shows a good level of generalization when applied to complex geometries, such as three-dimensional geometries reconstructed from medical images. \textcolor{black}{Finally, the model’s capability to perform agglomeration in heterogeneous domains is evaluated when integrated into a polytopal discontinuous Galerkin finite element solver.}
\end{abstract}

\section{Introduction}
Nowadays, the numerical solution of differential problems is central to all fields of engineering and applied sciences; moreover, it is expanding its boundaries to other fields that can benefit from numerical simulations. The applications in all these fields often involve complex physical domains, possibly featuring moving geometries and interfaces, heterogeneous media, immersed interfaces \cite{zonca_polygonal_2022}, and complex structures \cite{corti_numerical_2023,corti_discontinuous_2023}. The need to deal with such complex structures leads to the development of new and more sophisticated Finite Element Methods (FEM) that can employ general polygons and polyhedra (polytopal, in short) as grid elements for the numerical discretization of Partial Differential Equations (PDEs). We mention the Polyhedral Discontinuous Galerkin method \cite{antonietti_hp-version_2013,cangiani_hp-version_2014,antonietti_review_2016,cangianiHpVersionDiscontinuousGalerkin2017,antonietti_highorder_2021,antonietti_lymph_2024}, the Virtual Element Method \cite{beirao_da_veiga_basic_2013,beirao_da_veiga_hitchhikers_2014,beirao_da_veiga_virtual_2016,veiga_mixed_2016,antonietti_virtual_2022}, the Hybrid High-Order method \cite{di_pietro_arbitrary_2014,di_pietro_hybrid_2015,di_pietro_review_2016,di_pietro_hybrid_2020}, the mimetic finite difference method \cite{hyman_numerical_1997,brezzi_family_2005,brezzi_convergence_2005,da_veiga_mimetic_2014}, and the hybridizable discontinuous Galerkin method \cite{cockburn_superconvergent_2008,cockburn_superconvergent_2009,cockburn_unified_2009,cockburn_projection_2010}.
The flourishing of these methods takes with it the possibility of an intensive use of mesh agglomeration techniques that allow merging mesh elements to obtain coarser polytopal grids. However, the issue of efficiently performing agglomeration of polytopal meshes is still an open problem and object of intense research \cite{antonietti_agglomeration_2024,feder_r3mg_2024}. Mesh agglomeration is a powerful tool because it can reduce the computational complexity of the mesh while retaining a detailed geometrical description of the computational domain and/or interfaces. In scientific simulations, reducing the number of elements of a computational mesh can lead to faster simulations and quicker results while maintaining accuracy, for example, by locally increasing the polynomial order. Furthermore, agglomerated meshes form the basis for generating a hierarchy of (nested) coarser grids starting from a fine mesh of a complex physical domain of interest. Then, such grids' sequence, characterized by different refinement levels, can be employed in multi-grid solvers \cite{chan_convergence_1996,bassi_flexibilty_2012,bassi_agglomeration_2012,antonietti_multigrid_2015,antonietti_multigrid_2017,xu_algebraic_2017,antonietti_v_cycle_2019} to accelerate the numerical solution of the underlying algebraic system.\\

Because of the impossibility of the agglomerated grids to be used in the framework of classical FEMs, grid agglomeration is a relatively unexplored topic. Currently employed agglomeration methods rely on iterative algorithms based on graph partitioning of the connectivity matrix, which makes them not particularly feasible for complicated or heterogeneous domains or whenever the differential system models multiple physical processes that have to be appropriately taken into account in the agglomeration process. Moreover, during the agglomeration procedure, the quality preservation of the underlying mesh is fundamental because it might affect the overall stability and accuracy of the numerical method. Indeed, a suitably agglomerated mesh can give the same accuracy of the solution computed on the finer grid, but with fewer degrees of freedom. This leads to memory savings and lower computational loads.\\

In recent years the use of Machine Learning (ML) techniques to accelerate numerical methods has become widespread \cite{raissi_hidden_2018,hesthaven_non_intrusive_2018,ray_artificial_2018,raissi_physics_informed_2019,regazzoni_machine_2019,regazzoni_machine_2020,regazzoni_machine_2022,antonietti_accelerating_2023,antonietti_machine_2022,antonietti_refinement_2022}. In this work, we present a novel geometric deep-learning approach to address the problem of agglomerating three-dimensional meshes of heterogeneous domains. The proposed approach is based on Graph Neural Networks (GNNs). Indeed, we do not simply try to decide an \textit{a-priori} criterion to agglomerate the mesh, which would inevitably result in poor performance or high computational cost due to the impossibility of capturing all the possible configurations. Instead, our ML strategy automatically exploits and processes the huge amount of available data to learn only the distribution of the features of interest for the application, leading to high performance and computational efficiency. \textcolor{black}{Moreover, the strong generalization capabilities of GNN architectures in partitioning enable the agglomeration of large meshes using only small meshes in the training set, thereby significantly reducing training costs, as discussed in \cite{gatti_graph_2022}.}
\begin{figure}[t]
    \centering
    \includegraphics[width=.6\textwidth]{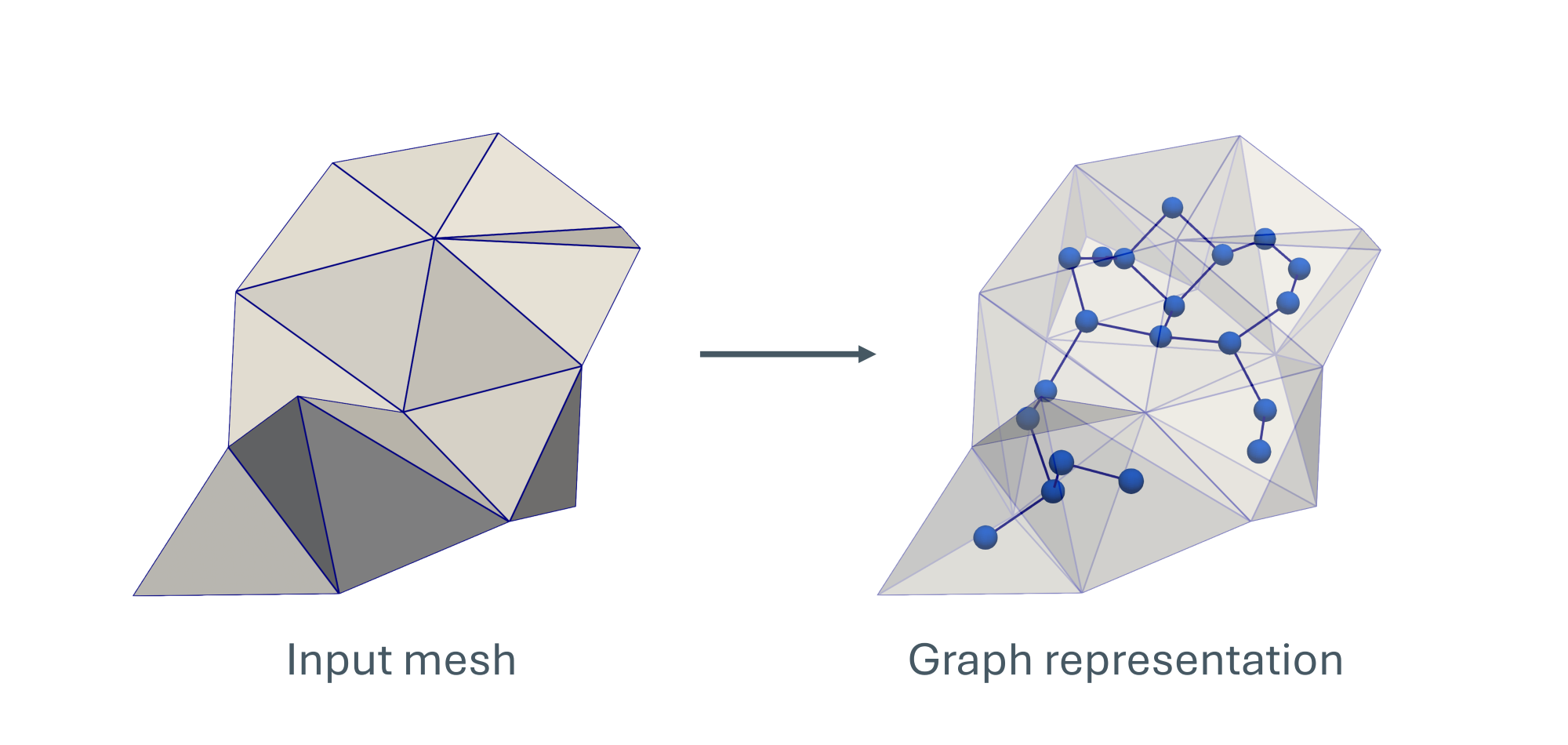}
    \caption{Graph extraction process: input mesh (left) and associated graph representation (right).}
    \label{fig:Graph_extraction}
\end{figure}
In our approach, the problem of mesh agglomeration is re-framed as a graph partitioning problem by exploiting the mesh structure. More precisely, the graph representation of the (dual) mesh is obtained by assigning a node to each mesh element. Edges connect pairs of nodes relative to adjacent mesh elements (see Figure~\ref{fig:Graph_extraction}). By exploiting such an equivalent representation, we can apply the proposed ML algorithm to solve a node classification problem, where each element is assigned to a cluster corresponding to a new polyhedral cell in the agglomerated mesh.
Furthermore, a powerful ability of the proposed GNN-based algorithm is the possibility to consider further physical properties of the computational domain during the process of mesh agglomeration. Indeed, by adding suitable input features and slightly modifying the network structure (particularly, changing the loss function), we can construct an algorithm that takes into account the heterogeneities of the underlying physical domain or possible physical features (e.g., embedded layers and/or microstructures). By exploiting this idea, we can develop a model that automatically agglomerates heterogeneous domains without the need to separately agglomerate regions with different material properties, with the need for a post-process to glue the obtained meshes to create the final one. Moreover, this new automatic algorithm can have an even bigger impact if the computational domain exhibits, for example, micro-structures and thus cannot be divided into simpler sub-regions. Further, we demonstrate that the proposed GNNs-based agglomeration algorithm has high online inference speed and can process naturally and simultaneously both the graph structure of the mesh and the geometrical information that can be associated with the graph nodes, such as the elements' areas/volumes or their centers' coordinates. The proposed GNNs-based agglomeration algorithm is then compared with the $k$-means \cite{macqueen_methods_1967} algorithm (another ML algorithm that can process only the geometrical information), and with METIS \cite{karypis_fast_1998}, a standard approach for graph partitioning that is meant to account only for the graph information on the mesh. The study extends the results of \cite{antonietti_agglomeration_2024}, which cannot be directly extended for agglomerating three-dimensional meshes or domains with complicated structures. Finally, we remark that the proposed algorithm is general, and the training can be done once and for all. More precisely, the models do not depend on the PDE under investigation or the numerical method adopted for the discretization.\\

The rest of the paper is organized as follows. In Section~\ref{sec:graphpart}, we present the GNN-agglomeration focusing on homogeneous domains in Section~\ref{sec:GNNhom}, and on heterogeneous ones in Section~\ref{sec:GNNhet}. 
In Section~\ref{sec:results_hom}, we present the numerical results focusing on homogeneous domains, first by measuring the effectiveness of the proposed strategies in terms of computational cost and quality metrics, then evaluating the generalization capabilities over a complex domain such as the human brain ventricles. Section~\ref{sec:results_het} focuses on the numerical results of heterogeneous domains by testing the quality of simple datasets and complex heterogeneous domains. \textcolor{black}{Section~\ref{sec:solver} presents numerical results demonstrating the performance of the proposed agglomeration strategy when the resulting meshes are used within PDE solvers.}
Finally, in Section~\ref{sec:conclusions}, we draw some conclusions.

\section{Mesh agglomeration algorithms via GNN}
\label{sec:graphpart}
\begin{figure}[t]
    \centering
    \includegraphics[width=\textwidth]{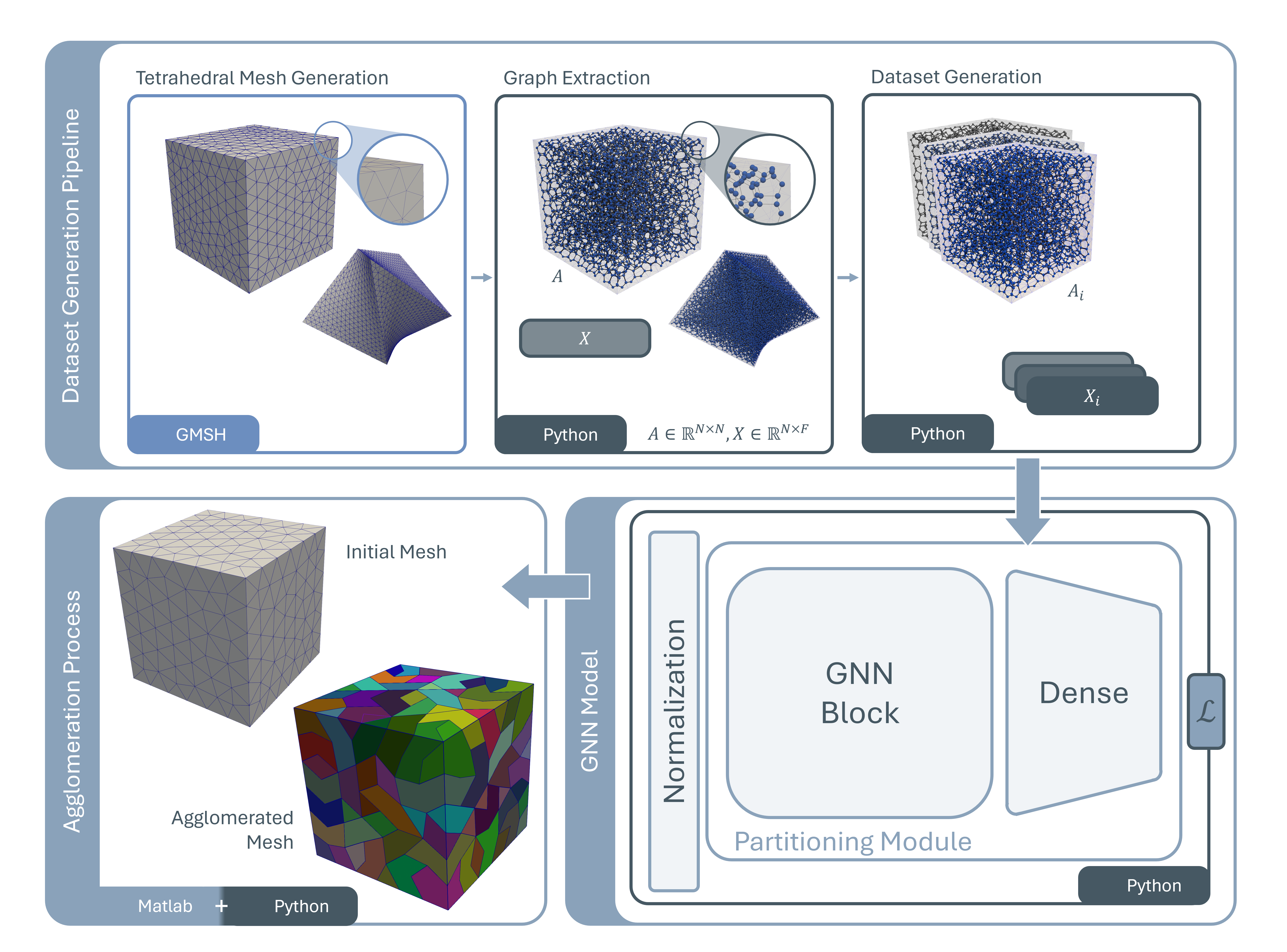}
    \caption{Schematic representation of the agglomeration pipeline.}
    \label{fig:pipeline}
\end{figure}

In this section, we introduce our agglomeration algorithm for three-dimensional meshes, discussing also the treatment of heterogeneous media. Most of the agglomeration algorithms tackle the problem by re-framing it into a graph partitioning problem by exploiting the connectivity structure of the mesh. In this work, we will focus on constructing an agglomeration algorithm based on recursive bisection of the associated graph. This strategy recursively bisects the input mesh's connectivity graph until the agglomerated elements have the desired size, as explained in Algorithm~\ref{alg:agglalgo}.
\begin{algorithm}[!hptb]
\caption{General mesh agglomeration strategy $\textsc{AGGLOMERATE}(\mathcal{T}_{h}, h^*)$}
\label{alg:agglalgo}
\begin{algorithmic}
\Require Mesh $\mathcal{T}_h$, target mesh size $h^*$, 
 and bisection model $\mathcal{M}$.
\Ensure Agglomerated mesh $\mathcal{T}_{h^*}$.
\If{$\diam(\mathcal{T}_h) \leq h^*$}
\Return $\mathcal{T}_h$
\Else
\State Extract the connectivity graph $G=(V,E)$ and features $X$ from $\mathcal{T}_h$.
\State Obtain the partition $Y$ by means of the bisection model $\mathcal{M}(G,X)$ (i.e., GNN, METIS or k-means).
\State Adjust partition $Y$.
\State Partition $\mathcal{T}_h$ into sub-meshes $\mathcal{T}_h^{(1)},\mathcal{T}_h^{(2)}$ according to $Y$.
\State $\mathcal{T}_{h^*}^{(1)} \leftarrow \textsc{agglomerate}(\mathcal{T}_h^{(1)},h^*)$
\State $\mathcal{T}_{h^*}^{(2)} \leftarrow \textsc{agglomerate}(\mathcal{T}_h^{(2)},h^*)$
\State$\mathcal{T}_{h^*} \leftarrow$ merge   $\mathcal{T}_{h^*}^{(1)}$, $\mathcal{T}_{h^*}^{(2)}$
\EndIf
\end{algorithmic}
\end{algorithm}
In the particular case of our strategy, we propose to employ GNN-based bisection models $\mathcal{M}$ that take as input a graph $G$ together with a set of features $X$ attached to each node (e.g. barycentric coordinates, area/volume and eventually physical parameters) and output the vector of probabilities $Y$ of each node to belong to cluster $1$ or $2$. This strategy has key advantages \textcolor{black}{compared to some of} the state-of-the-art methods like METIS \cite{karypis_fast_1998} and k-means \cite{macqueen_methods_1967}. More precisely:
\begin{itemize}
    \item It can process information about the connectivity graph of mesh elements, geometrical features of the elements, and physical properties of the modeled material. On the contrary, METIS is meant to process only graph-based information, and considering the geometrical and physical parameters would require modification of graph weights, which could not be straightforward to define in real applications. The k-means algorithm is typically used to process only geometrical information and does not take into account the connectivity matrix of the mesh.
    \item The online inference on new instances is much faster than METIS and k-means. Indeed, in our case, the main computational burden is faced during the training phase, while the other algorithms perform iterations online to provide an estimate of the solution. Moreover, the training phase is independent of the geometry and the underlying PDE model, so it is performed once and for all.
    \item  The GNN automatically processes additional information without requiring the user to know how that information can benefit performance. On the contrary, for METIS and k-means, this has to be explicitly modeled.
\end{itemize}
\textcolor{black}{We point out that in this study, we restrict our comparison to the k-means and METIS algorithms. A more comprehensive benchmarking (in terms of robustness and scalability) with other advanced, domain-specific partitioning strategies such as Coarsening Schemes for Graph Partitioning (e.g. \cite{Ron2011407, Safro2012369}), or SCOTCH \cite{pellegrini:hal-00770422,chevalier:hal-00402893}, which are known for their effectiveness in specific scenarios, falls beyond the scope of this work but represents an important objective for future research.}
Notice that GNN-driven agglomeration techniques were first proposed in \cite{antonietti_agglomeration_2024}, for a two-dimensional setting and underlying homogeneous domains; however, the approach in \cite{antonietti_agglomeration_2024} was limited to two-dimensional cases. The main novelty of the extension considered here is the possibility to handle also three-dimensional settings and properly take into account heterogeneities and/or embedded interfaces. More precisely, to handle the variability of three-dimensional meshes, which are in general much more complex, we propose suitable improvements with respect to the algorithm proposed in \cite{antonietti_agglomeration_2024}. Those improvements are both related to the pipeline and the GNN’s architecture. Moreover, when considering a three-dimensional setting, training and inference need to be feasible with much larger datasets, due to the broader amount of elements that three-dimensional meshes contain, even just in simple geometries. Moreover, we also propose a GNN architecture that can treat the agglomeration of heterogeneous domains, changing the loss function of the network. In Figure~\ref{fig:pipeline}, we report a schematic graphical conceptualization of the workflow of the proposed GNN agglomeration algorithm. In all cases, the dataset of meshes is created using GMSH \cite{geuzaine_gmsh_2009}. Then, the graph is extracted for each mesh, and the dataset is created by bundling graph adjacencies and feature matrices together. Next, the GNN model is trained in an unsupervised learning setting using a PyTorch Geometric \cite{paszke_pytorch_2019} implementation. Finally, the trained model is loaded in evaluation mode and used to perform the graph bisection iteratively, obtaining a clustered graph, which is then post-processed to build the agglomerated mesh.
\subsection{Architectures of GNN-agglomeration models for homogeneous domains}
\label{sec:GNNhom}
In this section, we present the proposed GNN architecture for mesh agglomeration for homogeneous domains. In particular, we describe the proposed architecture (GNN-AGGL-Enhanced), comparing it with the most straightforward possible three-dimensional extension of the one proposed in \cite{antonietti_agglomeration_2024} (GNN-AGGL-Base). Indeed, graph partitioning is a problem that focuses on the structure of the graph itself, regardless of any specific spatial representation. This implies that for graphs in two or three dimensions, the principles and algorithms for graph partitioning remain the same. At the same time, suitable modifications to handle three-dimensional data are needed. However, from the comparison of the result reported in Section~\ref{sec:results_hom} below, it is evident that there is a need for a more complex GNN architecture, as the GNN-AGGL-Enhanced one, for the agglomeration task in three dimensions.
\par
The proposed models (GNN-AGGL-Base and GNN-AGGL-Enhanced) are both based on two GNN SAGE-Base models, extending the ideas of \cite{antonietti_agglomeration_2024} to three-dimensional space. Both graph-bisection models perform a classification task, taking as input the graph $G=(V, E)$ and the features related to its nodes $X\in\mathbb{R}^{N\times4}$. The columns of the input feature matrix $X$ at the $i$-th row correspond to the information of the $i$-th mesh element. In particular, we have:
\begin{equation*}
    X = 
\begin{bmatrix}
x_1     & y_1       & z_1       & V_1      \\
\vdots  & \vdots    & \vdots    & \vdots   \\
x_i     & y_i       & z_i       & V_i      \\
\vdots  & \vdots    & \vdots    & \vdots   \\
x_N     & y_N       & z_N       & V_N 
\end{bmatrix},
\end{equation*}
where $(x_i,y_i,z_i)$ are the $(x,y,z)$-coordinates of the $i$-th element's centroid, while the $V_i$ is the element's volume. The output is $Y \in \mathbb{R}^{N\times2}$, where the first column represents the probability of belonging to the class $S_1$, while the second is the probability of being in the class $S_2$ of the partition. The characteristics of the resulting networks are reported in Table~\ref{tab:GNNs}. 
We underline that the number of parameters of the GNN-AGGL-Enhanced model is $\sim 110\mathrm{K}$, with respect to the $\sim 28\mathrm{K}$ of the GNN-AGGL-Base model. However, as will be discussed in Section~\ref{sec:results_hom}, this level of complexity is required to tame the inherited geometric complexity of the three-dimensional setting correctly.
\begin{table}[t]
    \centering
    \begin{tabular}{c|C C}
     & \textbf{GNN-AGGL-Base} & \textbf{GNN-AGGL-Enhanced}\\
     \hline & \\[-10pt]
     & Normalization of the features in $X$: & Normalization of the features in $X$: \\[2pt]  & \;\textbullet\; Coordinates rescaled in $[-1,1]$. & \;\textbullet\; Coordinates standardization: \\
     \textbf{Normalization} & & $\quad\quad$ Mean: $\mu = 0$ - Variance: $\sigma^2=1$ \\[2pt]
     \textbf{Layer} & \;\textbullet\; Volumes rescaled in $[0,1]$. & \;\textbullet\; Volumes rescaled in $[0,1]$. \\[2pt]
         & Rotation to align the maximum stretch direction with $x$-axis. & Rotation to align the maximum stretch direction with $x$-axis. \\[2pt]
     \hline & \\[-10pt]
     \textbf{SAGEConv} & \;\textbullet\; 4 layers of $64$ neurons. & \;\textbullet\; 4 layers of $128$ neurons. \\[2pt]
     \textbf{Layer} & \;\textbullet\; Activation function: $\tanh$. & \;\textbullet\; Activation function: $\tanh$. \\[2pt]
    \hline & \\[-10pt]
     \textbf{Linear Layer} & 3 decreasing-layers of $32-8-2$ neurons. & 4 decreasing-layers of $64-32-8-2$ neurons. \\[2pt]
     \hline & \\[-10pt]
     \multirow{2}{*}{\textbf{SoftMax Layer}} & \multicolumn{2}{>{\arraybackslash}p{14cm}}{Single layer with softmax function to bisect the result into the two classes $S_1$ and $S_2$, to obtain the resulting matrix $Y$.}\\[2pt]
     \hline
    \end{tabular}
    \caption{Main features of the two GNN architectures with a comparison of the specific features for each layer.}
    \label{tab:GNNs}
\end{table}
%
\subsubsection*{Training procedure}
\begin{figure}[t]
    \centering
    \includegraphics[width=\textwidth]{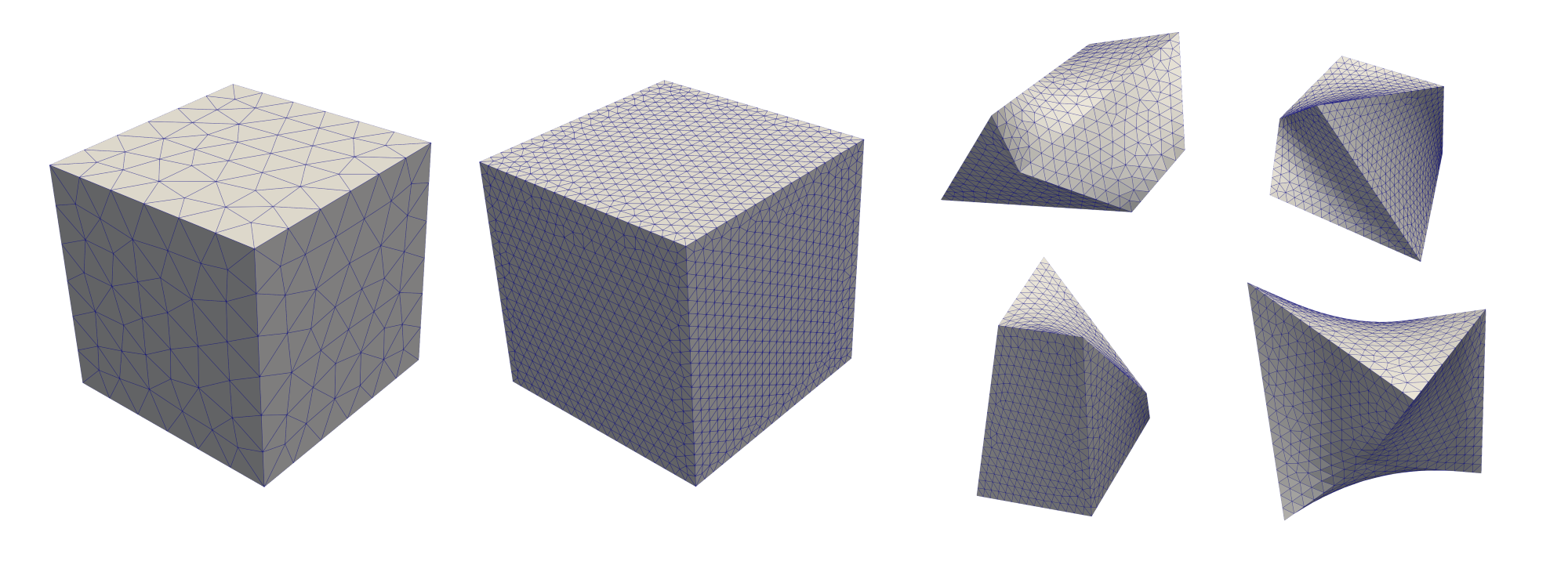}
    \caption{Examples of training dataset's meshes: meshes of a unitary cube (first two on the left) and a mesh of portions of the unitary cube (four ones on the right).}
    \label{fig:plot_mesh}
\end{figure}

The training of the models has been performed in an unsupervised learning setting. The loss function that is used is the \textit{expected normalized cut}, defined as:
\begin{equation}
    \mathcal{L}(Y,G) = \sum_{(i,j)\in E}\left(\dfrac{Y_{i1}(1-Y_{j1})}{\Gamma_1}+\dfrac{Y_{i2}(1-Y_{j2})}{\Gamma_2}\right),
\end{equation}
where $\Gamma_k = \sum_{i\in V} Y_{ik}D_i$ and $D$ is the vector of degrees of nodes. For details regarding the construction of the loss function, we refer to \cite{antonietti_agglomeration_2024}.
Concerning the construction of the dataset, we generate a training dataset of $500$ tetrahedral meshes ($125$ meshes of the unitary cube and $375$ meshes of suitable portions of them). All the meshes are generated using GMSH \cite{geuzaine_gmsh_2009}. We adopt a split between the training set and validation set, with the training set comprising $80\%$ and the validation set comprising $20\%$ of the total dataset. As a consequence, the training set is composed of $400$ grids ($100$ meshes of the unitary cube and $300$ meshes of portions). At the same time, the validation dataset is made by the remaining $100$ grids. To guarantee good agglomeration results with our GNN models, we use a dataset with a large variety of meshes, differentiating the number of mesh elements. In the training and validation datasets, we start from coarse meshes of $600$ elements, until around $42\,000$ elements of mesh in the finest one. In Figure~\ref{fig:plot_mesh}, we report some examples of meshes from the training set.
Due to the low dimension of such datasets, we also adopt the two strategies to avoid overfitting effects: 
\begin{itemize}
\item Regularization, summing a $L^2$-regularization term to the normalized cut loss, leading to the following loss definition:
\begin{equation}
\label{eq:loss_original}
 \hat{\mathcal{L}}(W,Y,G) = \mathcal{L}(Y,G) + \lambda \|W\|^2_2.
\end{equation}
\item Data-augmentation, by implementing a rotation of each sample by a random angle at each training iteration. This approach makes the model more robust with respect to the input mesh configurations.
\end{itemize}
The training process has been performed by adopting an Adam optimizer. All the details about the training hyperparameters for the two models are reported in Table~\ref{tab:training_hyper}.
\begin{table}
    \centering
    \begin{tabular}{|c|>{\centering\arraybackslash}p{4.5cm}|>{\centering\arraybackslash}p{4.5cm}|}
    \hline
     & \textbf{GNN-AGGL-Base} & \textbf{GNN-AGGL-Enhanced}\\
     \hline
     \textbf{Number of epochs} & $300$ Epochs & $400$ Epochs \\
     \hline 
     \textbf{Batch size} & $4$ Samples & $4$ Samples \\
     \hline 
     \textbf{Learning rate} ($\gamma$) & $10^{-5}$ & $10^{-4}$ \\
     \hline 
     \textbf{Weight-decay parameter} ($\lambda$) & $10^{-5}$ & $10^{-5}$ \\
     \hline
    \end{tabular}
    \caption{Training hyperparameters adopted for the two proposed GNN models.}
    \label{tab:training_hyper}
\end{table}
\begin{figure}[t]
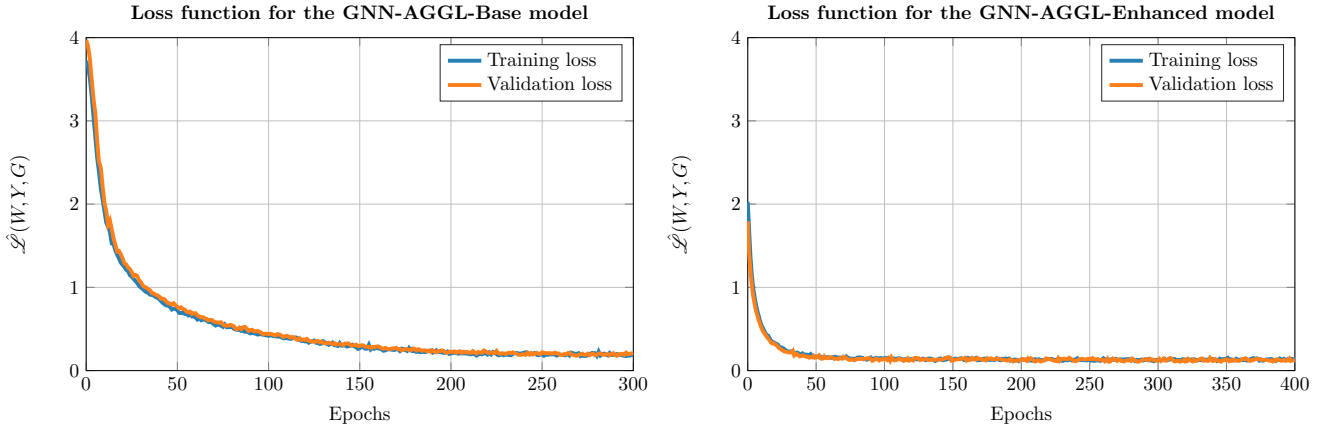

    \begin{subfigure}[b]{0.5\textwidth}
            \resizebox{\textwidth}{!}{\input{Training2DGNN.tex}}
    \end{subfigure}%
    \begin{subfigure}[b]{0.5\textwidth}
            \resizebox{\textwidth}{!}{\input{Training3DGNN.tex}}
    \end{subfigure}%
    \caption{Training and validation loss functions for the GNN-AGGL-Base model (left) and the GNN-AGGL-Enhanced model (right).}
    \label{fig:loss3D}
\end{figure}
In Figure~\ref{fig:loss3D}, we report the values of the loss functions for the training and validation sets for the GNN-AGGL-Base and the GNN-AGGL-Enhanced models. We can observe that using the GNN-AGGL-Enhanced model, the loss function reaches lower values on both datasets. Moreover, the training procedure is more stable compared to the two-dimensional setting \cite{antonietti_agglomeration_2024}. This is probably a consequence of the fact that the grids of the three-dimensional training dataset are made only of tetrahedra and not of general shape.
\subsection{Architectures of GNN-agglomeration model for heterogeneous domains}
\label{sec:GNNhet}
In this section, we further extend the proposed approach to perform mesh agglomeration on heterogeneous domains with different physical parameters. Specifically, this model can be exploited in different settings such as:
\begin{itemize}
\item Heterogeneous domains, where the main goal of the algorithm is to maintain as much as possible the distinction between the different material properties, still ensuring an acceptable shape for the agglomerated elements.
\item Domains with micro-structures, which is a more complex situation because there is no precise criterion to perform the agglomeration, and the algorithm aims to understand how to group the elements to obtain a better grid for numerical methods.
\end{itemize}
\textcolor{black}{We remark that graph partitioning algorithms like METIS or k-means are not able to perform agglomeration automatically, respecting the underlying heterogeneities and/or the embedded interfaces/microstructures.} The only possibility of performing this task with \textcolor{black}{the above-mentioned} methods is to do the agglomeration separately on the different regions and then merge them. However, the merging process of separate meshes is complex and time-consuming. Our model automates the procedure by exploiting the abilities of GNNs. \\

The GNN model proposed for the heterogeneous case is an extension of the GNN-AGGL-Enhanced model previously discussed. The most relevant change in the GNN architecture is the addition of a new input feature representing the physical parameter $\rho_i$ of the region to which the $i$-th element belongs. In this way we obtain an input feature matrix $X \in \mathbb{R}^{N\times5}$:
\begin{equation*}
    X = 
\begin{bmatrix}
x_1     & y_1       & z_1       & V_1      & \rho_1\\
\vdots  & \vdots    & \vdots    & \vdots   & \vdots\\
x_i     & y_i       & z_i       & V_i      & \rho_i\\
\vdots  & \vdots    & \vdots    & \vdots   & \vdots\\
x_N     & y_N       & z_N       & V_N      & \rho_N 
\end{bmatrix},
\end{equation*}
where $(x_i,y_i,z_i)$ are the $(x,y,z)$-coordinates of the $i$-th element's center, while the $V_i$ is the element's volume. As explained for the other quantities, the parameter $\rho_i$ would be re-scaled to obtain values between $0$ and $1$. This practice reduces the variability of the data, making it easier for the network to learn the relative importance of this parameter on the final result.
It is well known that neural networks often struggle to handle sharp changes in input data that may lead to gradient instability and hinder the model's generalization ability. Adding a new pre-processing operation can easily overcome this obstacle for heterogeneous domains. For this reason, we substitute in the normalization layer the last input feature with the average of such values in the adjacent elements. This step allows for avoiding sharp changes at the boundaries of distinct physical regions by smoothing out this transition.
The output of the GNN remains unchanged. It is $Y \in \mathbb{R}^{N\times2}$, where the first column represents the probability of belonging to $S_1$, while the second is the probability of being in $S_2$. 

The idea behind the agglomeration algorithm is similar to that described in Algorithm~\ref{alg:agglalgo}. However, we introduce some modifications to better deal with the specific task of this model. Specifically, we add an additional check before adding the agglomerated element to the final mesh to ensure the polytopal element connectivity. If the resulting agglomerate is not connected (an efficient algorithm for finding the connected components of a graph is reported in \cite{tarjan_depth_1972}), we perform an extra agglomeration step. The need to introduce this check is due to the network's tendency to keep elements with the same physical parameter together, even when they are not connected.
\subsubsection*{Training procedure}
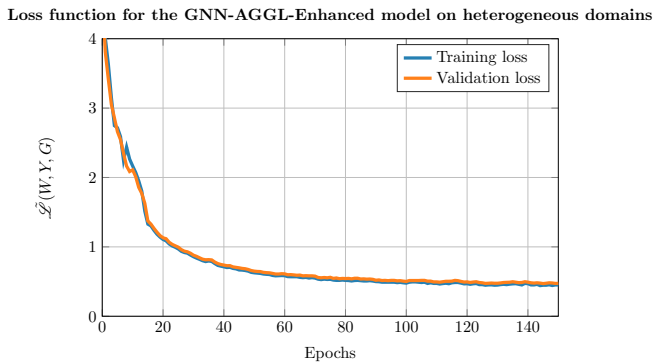
\begin{figure}[t]
\centering
\resizebox{0.5\textwidth}{!}{\definecolor{Validation}{rgb}{1.00000,0.50000,0.10000}%
\definecolor{Training}{rgb}{0.15000,0.50000,0.70000}%
\begin{tikzpicture}

\begin{axis}[%
width=3.875in,
height=2.36in,
at={(2.6in,1.099in)},
scale only axis,
xmin=0,
xmax=150,
xminorticks=true,
xlabel = {Epochs},
ylabel = {$\tilde{\mathcal{L}}(W,Y,G)$},
ymin=0,
ymax=4,
yminorticks=true,
axis background/.style={fill=white},
title style={font=\bfseries},
title={Loss function for the GNN-AGGL-Enhanced model on heterogeneous domains},
xmajorgrids,
xminorgrids,
ymajorgrids,
yminorgrids,
legend style={legend cell align=left, align=left, draw=white!15!black}
]
                
\addplot [color=Training, line width=2.0pt]
  table[row sep=crcr]{%
0.00000 4.88509 \\
1.00000 4.01271 \\
2.00000 3.63517 \\
3.00000 3.15308 \\
4.00000 2.74719 \\
5.00000 2.71674 \\
6.00000 2.58996 \\
7.00000 2.24676 \\
8.00000 2.43857 \\
9.00000 2.26558 \\
10.00000 2.16562 \\
11.00000 2.07253 \\
12.00000 1.95081 \\
13.00000 1.80467 \\
14.00000 1.51351 \\
15.00000 1.32827 \\
16.00000 1.30264 \\
17.00000 1.23871 \\
18.00000 1.18308 \\
19.00000 1.14081 \\
20.00000 1.10684 \\
21.00000 1.08709 \\
22.00000 1.04272 \\
23.00000 1.01399 \\
24.00000 0.99428 \\
25.00000 0.97029 \\
26.00000 0.93614 \\
27.00000 0.91859 \\
28.00000 0.90934 \\
29.00000 0.88273 \\
30.00000 0.85829 \\
31.00000 0.83818 \\
32.00000 0.82046 \\
33.00000 0.80063 \\
34.00000 0.78793 \\
35.00000 0.79659 \\
36.00000 0.79027 \\
37.00000 0.75970 \\
38.00000 0.73309 \\
39.00000 0.72216 \\
40.00000 0.71266 \\
41.00000 0.70237 \\
42.00000 0.70519 \\
43.00000 0.68952 \\
44.00000 0.68293 \\
45.00000 0.66739 \\
46.00000 0.66539 \\
47.00000 0.65955 \\
48.00000 0.65042 \\
49.00000 0.63271 \\
50.00000 0.62619 \\
51.00000 0.62200 \\
52.00000 0.61770 \\
53.00000 0.60872 \\
54.00000 0.60436 \\
55.00000 0.60184 \\
56.00000 0.58854 \\
57.00000 0.58253 \\
58.00000 0.58472 \\
59.00000 0.58691 \\
60.00000 0.58469 \\
61.00000 0.57144 \\
62.00000 0.56952 \\
63.00000 0.56945 \\
64.00000 0.56967 \\
65.00000 0.56073 \\
66.00000 0.55636 \\
67.00000 0.56021 \\
68.00000 0.55074 \\
69.00000 0.55636 \\
70.00000 0.54997 \\
71.00000 0.53457 \\
72.00000 0.52802 \\
73.00000 0.53497 \\
74.00000 0.53018 \\
75.00000 0.53255 \\
76.00000 0.52146 \\
77.00000 0.51709 \\
78.00000 0.51709 \\
79.00000 0.52146 \\
80.00000 0.51709 \\
81.00000 0.51502 \\
82.00000 0.50850 \\
83.00000 0.51488 \\
84.00000 0.51709 \\
85.00000 0.51271 \\
86.00000 0.50836 \\
87.00000 0.50836 \\
88.00000 0.51273 \\
89.00000 0.50752 \\
90.00000 0.49942 \\
91.00000 0.49755 \\
92.00000 0.49091 \\
93.00000 0.49091 \\
94.00000 0.49091 \\
95.00000 0.49091 \\
96.00000 0.49093 \\
97.00000 0.48442 \\
98.00000 0.48865 \\
99.00000 0.48218 \\
100.00000 0.47762 \\
101.00000 0.49087 \\
102.00000 0.49527 \\
103.00000 0.49526 \\
104.00000 0.49091 \\
105.00000 0.48637 \\
106.00000 0.49107 \\
107.00000 0.47306 \\
108.00000 0.48219 \\
109.00000 0.47345 \\
110.00000 0.47354 \\
111.00000 0.46473 \\
112.00000 0.47122 \\
113.00000 0.47580 \\
114.00000 0.46909 \\
115.00000 0.48219 \\
116.00000 0.48655 \\
117.00000 0.48677 \\
118.00000 0.47770 \\
119.00000 0.46686 \\
120.00000 0.46908 \\
121.00000 0.46069 \\
122.00000 0.46783 \\
123.00000 0.47193 \\
124.00000 0.45594 \\
125.00000 0.45598 \\
126.00000 0.44728 \\
127.00000 0.45164 \\
128.00000 0.45256 \\
129.00000 0.44727 \\
130.00000 0.44936 \\
131.00000 0.45614 \\
132.00000 0.46044 \\
133.00000 0.45600 \\
134.00000 0.46247 \\
135.00000 0.46697 \\
136.00000 0.46911 \\
137.00000 0.46043 \\
138.00000 0.45148 \\
139.00000 0.47366 \\
140.00000 0.46945 \\
141.00000 0.45170 \\
142.00000 0.45159 \\
143.00000 0.45969 \\
144.00000 0.44161 \\
145.00000 0.44727 \\
146.00000 0.44728 \\
147.00000 0.45385 \\
148.00000 0.44317 \\
149.00000 0.45145 \\
150.00000 0.45145 \\
};
\addlegendentry{Training loss}

\addplot [color=Validation, line width=2.0pt]
  table[row sep=crcr]{%
0.00000 4.69309 \\
1.00000 3.81105 \\
2.00000 3.41787 \\
3.00000 3.05414 \\
4.00000 2.85829 \\
5.00000 2.66028 \\
6.00000 2.54565 \\
7.00000 2.34808 \\
8.00000 2.16466 \\
9.00000 2.08586 \\
10.00000 2.11021 \\
11.00000 2.01137 \\
12.00000 1.85708 \\
13.00000 1.77010 \\
14.00000 1.62718 \\
15.00000 1.37296 \\
16.00000 1.32868 \\
17.00000 1.26707 \\
18.00000 1.21174 \\
19.00000 1.16268 \\
20.00000 1.13071 \\
21.00000 1.11740 \\
22.00000 1.06665 \\
23.00000 1.03642 \\
24.00000 1.01563 \\
25.00000 0.99595 \\
26.00000 0.96066 \\
27.00000 0.93818 \\
28.00000 0.92912 \\
29.00000 0.90880 \\
30.00000 0.88009 \\
31.00000 0.86378 \\
32.00000 0.84250 \\
33.00000 0.82245 \\
34.00000 0.81406 \\
35.00000 0.81600 \\
36.00000 0.81209 \\
37.00000 0.78171 \\
38.00000 0.75928 \\
39.00000 0.74622 \\
40.00000 0.73737 \\
41.00000 0.72873 \\
42.00000 0.72919 \\
43.00000 0.71150 \\
44.00000 0.70472 \\
45.00000 0.69809 \\
46.00000 0.68966 \\
47.00000 0.68502 \\
48.00000 0.67443 \\
49.00000 0.65889 \\
50.00000 0.65024 \\
51.00000 0.64582 \\
52.00000 0.64389 \\
53.00000 0.63272 \\
54.00000 0.62836 \\
55.00000 0.62597 \\
56.00000 0.61775 \\
57.00000 0.60879 \\
58.00000 0.60849 \\
59.00000 0.61330 \\
60.00000 0.61086 \\
61.00000 0.59788 \\
62.00000 0.60021 \\
63.00000 0.59587 \\
64.00000 0.59079 \\
65.00000 0.59344 \\
66.00000 0.58239 \\
67.00000 0.58469 \\
68.00000 0.58255 \\
69.00000 0.58255 \\
70.00000 0.57828 \\
71.00000 0.55846 \\
72.00000 0.55419 \\
73.00000 0.55854 \\
74.00000 0.55419 \\
75.00000 0.56110 \\
76.00000 0.54545 \\
77.00000 0.55169 \\
78.00000 0.54115 \\
79.00000 0.54545 \\
80.00000 0.54326 \\
81.00000 0.54564 \\
82.00000 0.53910 \\
83.00000 0.54763 \\
84.00000 0.54786 \\
85.00000 0.53667 \\
86.00000 0.53891 \\
87.00000 0.53454 \\
88.00000 0.53673 \\
89.00000 0.53215 \\
90.00000 0.51906 \\
91.00000 0.52357 \\
92.00000 0.51490 \\
93.00000 0.51709 \\
94.00000 0.51277 \\
95.00000 0.51709 \\
96.00000 0.51703 \\
97.00000 0.51055 \\
98.00000 0.51267 \\
99.00000 0.50619 \\
100.00000 0.50168 \\
101.00000 0.51492 \\
102.00000 0.51491 \\
103.00000 0.51491 \\
104.00000 0.51491 \\
105.00000 0.51709 \\
106.00000 0.51705 \\
107.00000 0.50624 \\
108.00000 0.50619 \\
109.00000 0.49519 \\
110.00000 0.49524 \\
111.00000 0.49091 \\
112.00000 0.49758 \\
113.00000 0.50395 \\
114.00000 0.50400 \\
115.00000 0.50620 \\
116.00000 0.51702 \\
117.00000 0.51277 \\
118.00000 0.50222 \\
119.00000 0.49309 \\
120.00000 0.49531 \\
121.00000 0.48685 \\
122.00000 0.49326 \\
123.00000 0.49734 \\
124.00000 0.48647 \\
125.00000 0.48217 \\
126.00000 0.47334 \\
127.00000 0.47564 \\
128.00000 0.47875 \\
129.00000 0.47564 \\
130.00000 0.47563 \\
131.00000 0.48001 \\
132.00000 0.48655 \\
133.00000 0.48873 \\
134.00000 0.49040 \\
135.00000 0.49742 \\
136.00000 0.49091 \\
137.00000 0.48442 \\
138.00000 0.48655 \\
139.00000 0.49765 \\
140.00000 0.49122 \\
141.00000 0.48218 \\
142.00000 0.47804 \\
143.00000 0.48214 \\
144.00000 0.47804 \\
145.00000 0.47090 \\
146.00000 0.47573 \\
147.00000 0.48217 \\
148.00000 0.47955 \\
149.00000 0.47344 \\
150.00000 0.47344 \\
};
\addlegendentry{Validation loss}

\end{axis}
\end{tikzpicture}}
\caption{Training and validation loss functions for the GNN-AGGL-Enhanced model 
 for heterogeneous domains.}
\label{fig:loss3Dheter}
\end{figure}

Since this work is intended to be the first approach to GNN methods for heterogeneous mesh agglomeration, we use a training dataset composed only of heterogeneous domains with a small number of parameters and regions as a starting point. However, it is important to emphasize that this choice also allows us to perform agglomeration in the most general situations of computational domains with $S$ parameters or domains with micro-structures, as shown in the forthcoming numerical tests. Moreover, this choice allows us to evaluate the generalization capabilities of the model.
The training dataset contains $200$ tetrahedral meshes ($100$ with a single physical parameter and $100$ with two heterogeneous parameters). The validation dataset contains $50$ tetrahedral meshes ($25$ with a single physical parameter and $25$ with two different parameter regions). The choice of the number of elements in the dataset preserves the $80/20$ split between training and validation sets. The elements of the meshes in the two considered datasets are approximately $200$ for the coarsest mesh and $17\,000$ for the finest one.
The training of the models has been performed in an unsupervised learning setting. We slightly modify the loss function by adding a further term to consider the physical feature. Hence, the new loss function includes a term involving the physical parameters. Consider the matrix $P \in \mathbb{R}^{N\times2}$, which contains in the first column at the $i$-th row the "physical" parameter $p_i$ and in the second column its complement $1-p_i$ associated with the $i$-th mesh element. Then, the physical penalty part of the loss function reads:
\begin{equation}
    \label{eq:loss_penalty}
    \mathcal{P}(P,Y) = \sum P \odot Y,
\end{equation}
where $\odot$ denotes the element-wise multiplication and the summation is over all the entries of the resulting matrix. Finally, the loss function of the GNN is composed of the sum of the terms in Equations \eqref{eq:loss_original} and \eqref{eq:loss_penalty}:
\begin{equation}
    \tilde{\mathcal{L}}(W,P,Y,G) = \alpha\mathcal{P}(P,Y) + \beta\left(\mathcal{L}(W,Y,G) + \lambda\|W\|_2^2\right).
\end{equation}
The proposed loss function contains two parameters, $\alpha$ and $\beta$, to correctly balance the relevance of the normalized cut term and the term related to physical parameters. A hyperparameter tuning procedure has been used to balance the two terms. Indeed, by reducing the normalized cut's importance, the model separates the physical classes correctly, but we obtain low-quality agglomerated elements. Instead, by reducing the significance of the physical parameter term, we fail to account for the physical classes. The optimal values we found for that model are $\alpha = 1.28$ and $\beta = 2.2 \times 10^{-4}$.
The training has been performed by adopting $150$ epochs of the Adam optimizer with a batch size of $4$. The learning rate is set to be $\gamma=10^{-4}$, and the weight decay parameter is $\lambda=10^{-5}$ and batch size $4$. In Figure~\ref{fig:loss3Dheter} we report the decay of the loss function $\tilde{\mathcal{L}}$ evaluated on training and validation datasets. \textcolor{black}{The GNN training was carried out on a workstation equipped with an Intel Core i7-1270P CPU (12 cores, 16 threads, base clock 2.2GHz) and 32GB RAM.}
\textcolor{black}{Concerning the training times, we have to consider a cost of $\simeq 500$ minutes for the complete training of the GNN-AGGL-Enhanced network. However, this cost is offline, so we need to perform it only once. This cost can be compensated with the lower online costs with respect to the METIS or k-means algorithms, as discussed in the next sections. Moreover, the training set can be limited, containing small meshes, due to the high generalization properties over larger meshes of the GNN architecture in partitioning \cite{gatti_graph_2022}. This is a large advantage of the proposed method, which can be applied to agglomerate large meshes without an increase in the training costs.}
\section{Numerical results: mesh agglomeration for homogeneous domains}
\label{sec:results_hom}
In this section, after the introduction of some quality metrics for polyhedral meshes, we introduce a test dataset to compare the performance of the GNN model with \textcolor{black}{graph-partitioning methods like METIS and k-means.}
Moreover, we present an application to a mesh of human brain ventricles. \\

First, let us define the quality metrics we use in the next test cases. These quality metrics are employed also in \cite{attene_benchmarking_2021,antonietti_agglomeration_2024}, but we adapt them to the three-dimensional setting:
\begin{itemize}
\item \textbf{Circle Ratio (CR)}: ratio between the radius of the largest contained sphere and the radius of the smallest sphere containing the polyhedron $K$:
\begin{equation}
\operatorname{CR}(K)=\dfrac{\underset{B(r) \subset K}{\max}\{r\}}{\underset{K \subset B(r)}{\min}\{r\}},
\end{equation}
where $B(r)$ is a ball of radius $r$. We point out that inscribed spheres tangent to each of the polyhedron's faces and circumscribed ones touching each of the polyhedron's vertices might not exist for non-convex polyhedra. For this reason, we instead use the largest sphere contained in the polyhedron as the inscribed one and the smallest sphere containing the polyhedron as the circumscribed one \cite{attene_benchmarking_2021}.
\item \textbf{Uniformity Factor (UF)}: the ratio between the diameter of an element $K$ and the mesh size $h$:
\begin{equation}
    \operatorname{UF}(K)=\frac{\diam(K)}{h},
\end{equation}
where $h=\max_{K \in \mathcal{T}_h}{\{\diam(K)\}}$.
\item \textbf{Volumes Difference (VD)}:
    \begin{equation}
        \operatorname{VD}(K)=\frac{|V(K)-V_\mathrm{target}|}{V_\mathrm{target}},
    \end{equation}
    where $V_\mathrm{target} = \frac{\sum_{K\in\mathcal{T}_h}V(K)}{N_\mathrm{el}}$, where $N_\mathrm{el}$ is the number of the agglomerated elements of the mesh. Then, the target volume is the volume of an element in the agglomerated mesh if all the portions have the same volume.
\end{itemize}
The CR and UF quality metrics assume values between $0$ and $1$. Moreover, if the elements are regular, these quality metrics are close to $1$. Finally, we introduce a metric to evaluate the performance of the agglomeration algorithm in terms of the difference between the volumes of the agglomerated mesh elements. VD assumes positive values, and if the mesh is quasi-uniform, these values are close to $0$.

\subsection{Test case 1: agglomeration capabilities assessment of a test dataset}
\label{sec:test_case_1}
We construct a test dataset of $100$ tetrahedral meshes with different levels of refinements, ranging from  $563$ elements for the coarsest mesh and $6\,404$ for the finest one. The dataset is composed as follows: 1) $30$ meshes of the unitary cube; 2) $70$ meshes of suitable portions of the unitary cube.
The computational domains are analogous to the ones used for training, and we refer to Figure~\ref{fig:plot_mesh} that reports a mesh for each type. The results in the following sections are obtained using this test dataset. The considered agglomerated meshes have approximately $128$ elements, corresponding to $7$ bisection steps.
\begin{figure}[t]
    \centering
    \includegraphics[width=0.7\textwidth]{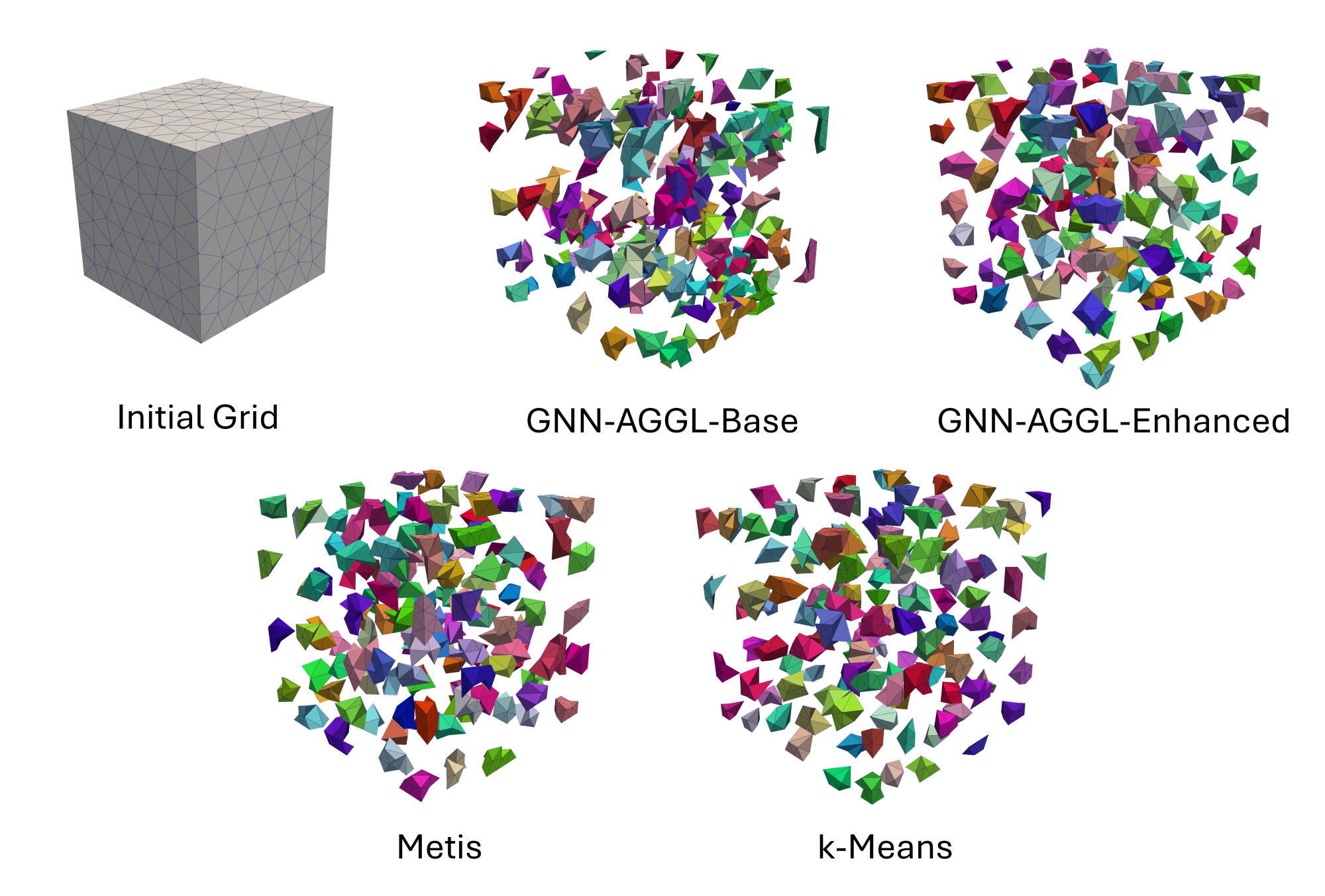}
    \caption{Test case 1: agglomerated grids obtained based on employing the GNN-AGGL-Base, the GNN-AGGL-Enhanced, METIS, and k-means algorithm.}    \label{fig:testcase1agglgrids}
\end{figure}
\begin{table}[t]
\centering
\begin{tabular}{|c|Q|Q|Q|}
\hline
\textbf{Agglomeration} & \textbf{Circle} & \textbf{Uniformity} & \textbf{Volume} \\
\textbf{Method} & \textbf{Ratio} & \textbf{Factor} & \textbf{Difference} \\ \hline
\textbf{GNN-AGGL-Base} 
& \cellcolor[RGB]{255, 68,0} 0.1342 & \cellcolor[RGB]{255,201,0} 0.4053 & \cellcolor[RGB]{255,144,0} 0.7166 \\ \hline
\textbf{GNN-AGGL-Enhanced} 
& \cellcolor[RGB]{255, 99,0} 0.1946 & \cellcolor[RGB]{118,255,0} 0.7676 & \cellcolor[RGB]{120,255,0} 0.2355 \\ \hline
\textbf{METIS}
& \cellcolor[RGB]{255, 92,0} 0.1799 & \cellcolor[RGB]{143,255,0} 0.7183 & \cellcolor[RGB]{179,255,0} 0.3517 \\ \hline
\textbf{k-Means} 
& \cellcolor[RGB]{255,136,0} 0.2662 & \cellcolor[RGB]{103,255,0} 0.7982 & \cellcolor[RGB]{167,255,0} 0.3270 \\ \hline
\end{tabular}
\caption{Test case 1: average values of the computed quality metrics (CR, UF, and VD) for the agglomerated grids reported in Figure~\ref{fig:testcase1agglgrids}, obtained with different agglomeration strategies (GNN-AGGL-Base, GNN-AGGL-Enhanced, METIS, k-Means). The colors are scaled starting from the minimum value of the metric (red) to the maximum one (green).}
\label{table:testcase1quality}
\end{table}
\begin{figure}[t!]
    \begin{subfigure}[b]{0.3\textwidth}
    \resizebox{\textwidth}{!}{\includegraphics{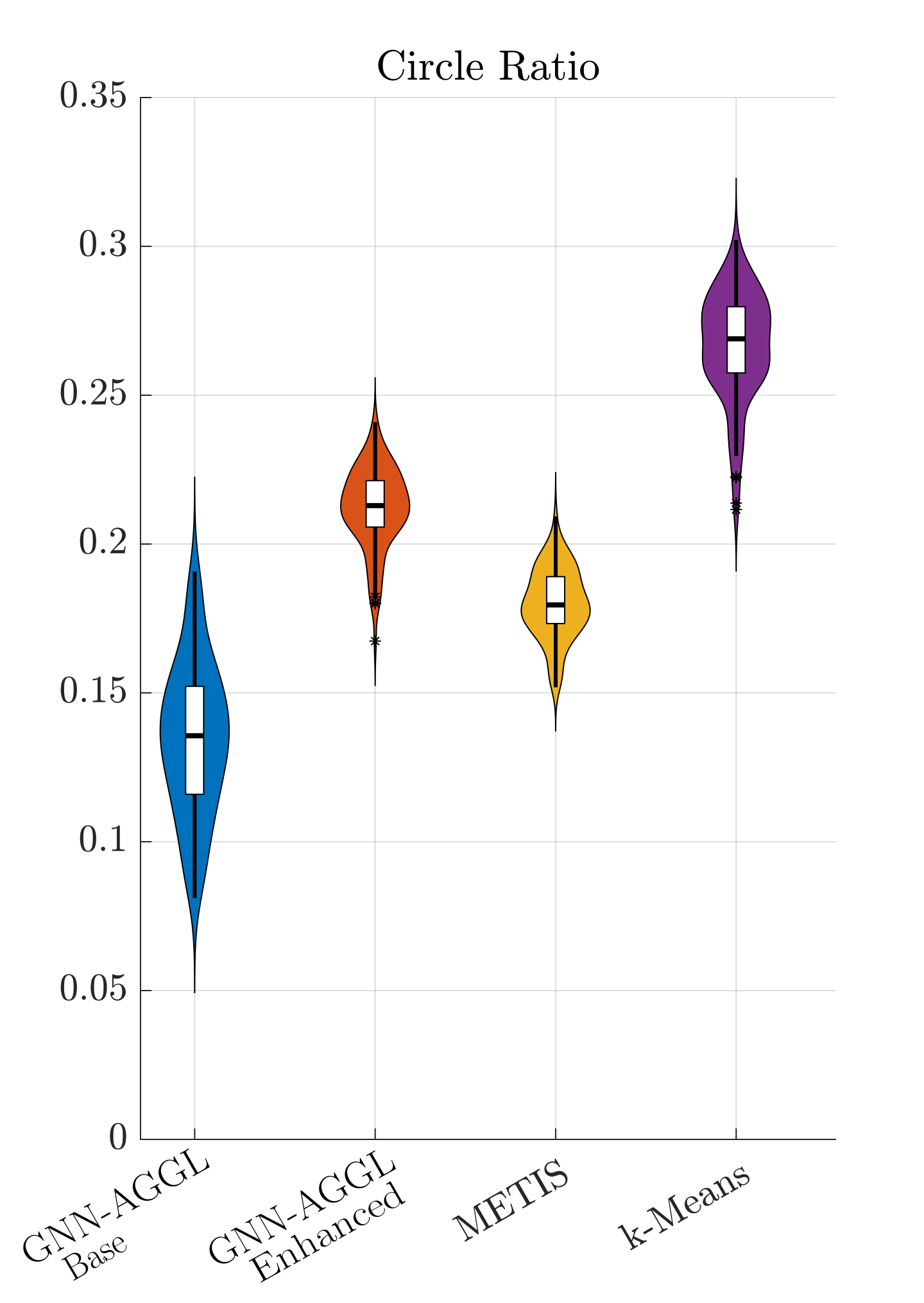}}
    \end{subfigure}
    \begin{subfigure}[b]{0.3\textwidth}
    \resizebox{\textwidth}{!}{\includegraphics{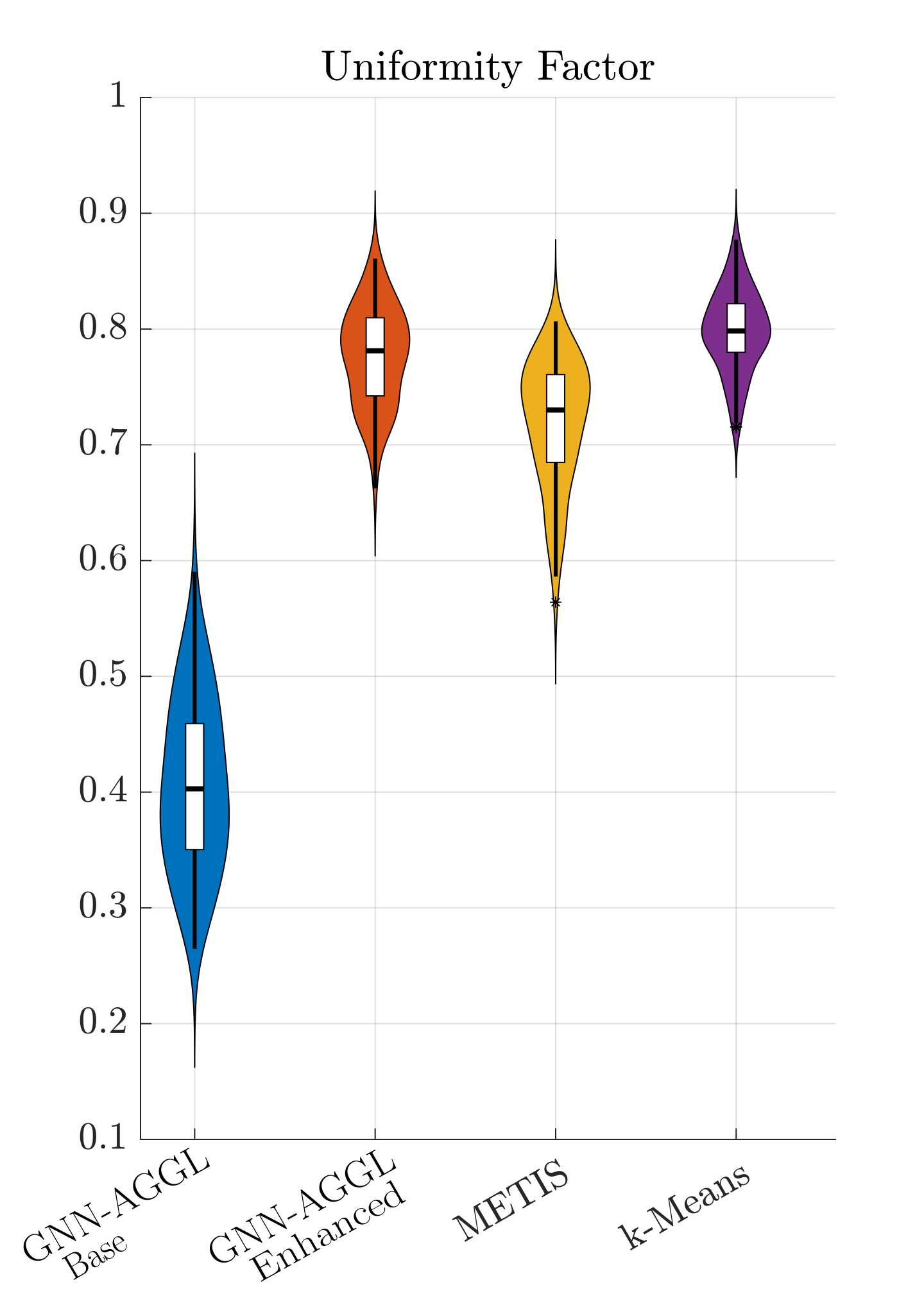}}
    \end{subfigure} 
    \begin{subfigure}[b]{0.3\textwidth}
    \resizebox{\textwidth}{!}{\includegraphics{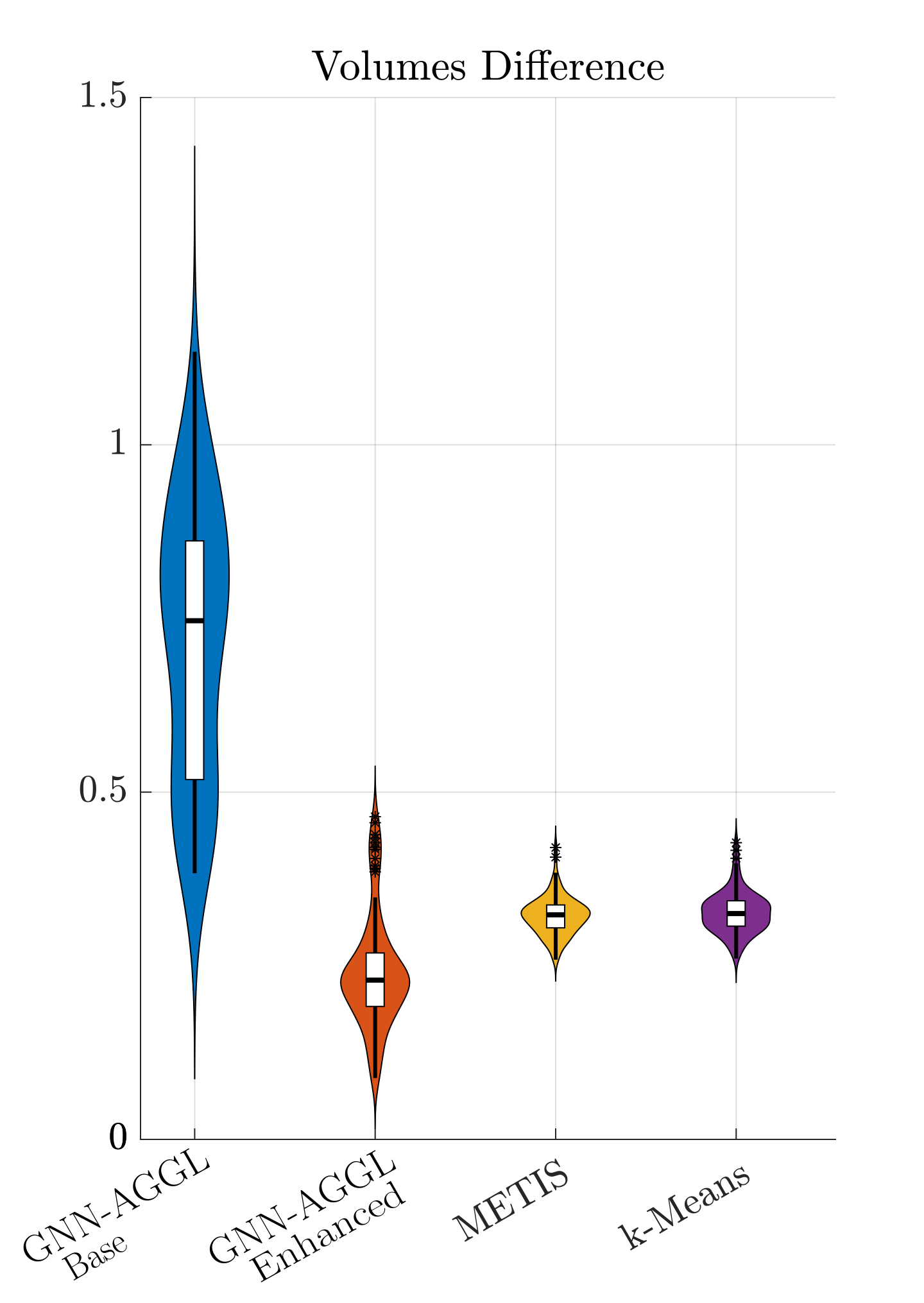}}
    \end{subfigure}%
    \caption{Test case 1: violin plots of the computed quality metrics (CR on the left, UF in the middle, and VD on the right) for the agglomerated grids reported in Figures \ref{fig:testcase1agglgrids}, obtained with different agglomeration strategies (GNN-AGGL-Base, GNN-AGGL-Enhanced, METIS, k-Means).}
    \label{fig:testcase1violin}
\end{figure}
In Figure~\ref{fig:testcase1agglgrids}, we report a mesh from the test dataset together with the agglomerated meshes obtained using four different partitioning methods, namely GNN-AGGL-Base, GNN-AGGL-Enhanced, METIS, and k-means. We plot the agglomerated meshes with an "exploded" view for visualization purposes. This allows us to see how tetrahedra are agglomerated together into the new mesh elements. We can see that GNN-AGGL-Enhanced and k-means grids look similar. On the contrary, METIS and GNN-AGGL-Base ones are of lower quality, due to the presence of stretched elements.
Table~\ref{table:testcase1quality} collects the average values over the full test set of the quality metrics (CR, UF, and VD) for the agglomerated grids obtained employing different agglomeration strategies (GNN-AGGL-Base, GNN-AGGL-Enhanced, METIS, and k-Means). Moreover, in Figure~\ref{fig:testcase1violin}, we report the violin plots relative to the distributions of the quality metrics for the considered agglomerated mesh. Looking at the results, it is clear that the GNN-AGGL-Base model is not able to create an agglomerated mesh of sufficiently good quality. As expected, due to the much greater variability in the three-dimensional context, the performance of this model is significantly lower than that obtained in the two-dimensional case \cite{antonietti_agglomeration_2024}. For this reason, we focus only on the GNN-AGGL-Enhanced model in the following test cases.
On the contrary, we can see that the GNN-AGGL-Enhanced model performs better than METIS. This can be easily explained because the latter considers only the information coming from the graph topology, while GNN-AGGL-Enhanced takes in input some geometrical properties of the mesh.
Observing also k-means, we can infer that the developed algorithm seems to be the best in terms of UF, while it is slightly worse than k-means in terms of shape regularity (CR index). However, we observe that the CR index provides low values for all the considered methods.\\
\textcolor{black}{We conclude this section by presenting some results on runtime performance.}
\begin{figure}[t]
    \centering
    \includegraphics[width=0.6\textwidth]{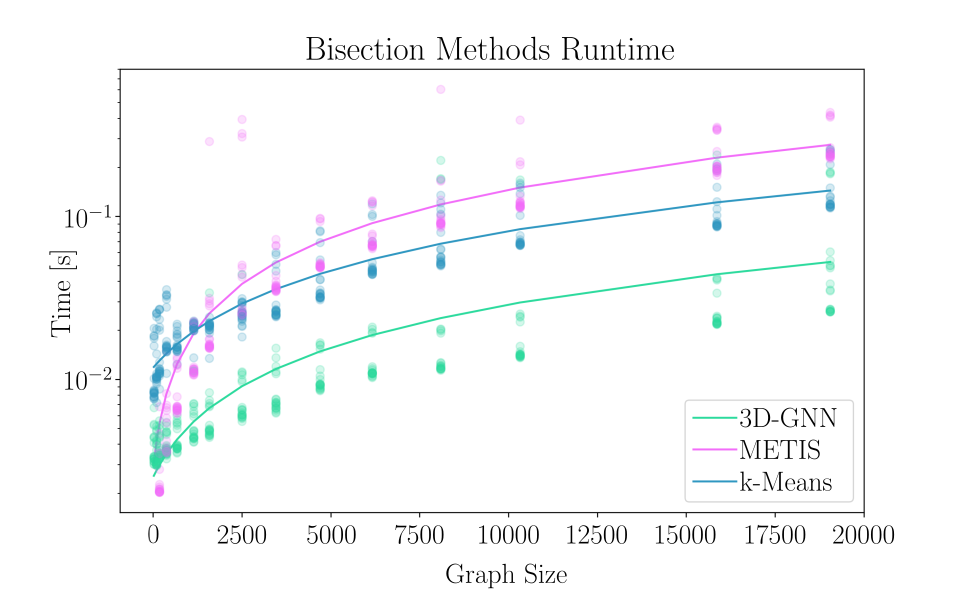}
    \caption{Test case 1: runtime performance for different graph bisection models (METIS, k-Means, and GNN-AGGL-Enhanced) as a function of the number of nodes in the connectivity graph of tetrahedral meshes.}
    \label{fig:runtimeplot}
\end{figure}
In Figure~\ref{fig:runtimeplot}, we report the runtime performance of the proposed GNN-AGGL-Enhanced algorithm, and compare it with the analogous quantities obtained with METIS and k-Means algorithms. On the abscissa, we show the graph size, i.e., the number of mesh elements, and on the ordinate, we plot the time required for the graph's bisection. The main advantage of our model with respect to the others is its non-iterative nature in performing graph bisection. The time benchmarking test set contains $15$ tetrahedral meshes of increasing size, from $24$ to $19 \ 049$ elements.
Due to the advantage mentioned above, the GNN-AGGL-Enhanced model seems to outperform both METIS and k-means in terms of runtime performance.

\subsection{Test case 2: agglomeration of a mesh of human brain ventricles}
\label{sec:agglomeration_ventricles}
\begin{figure}[t]
    \centering
  \includegraphics[width=\textwidth]{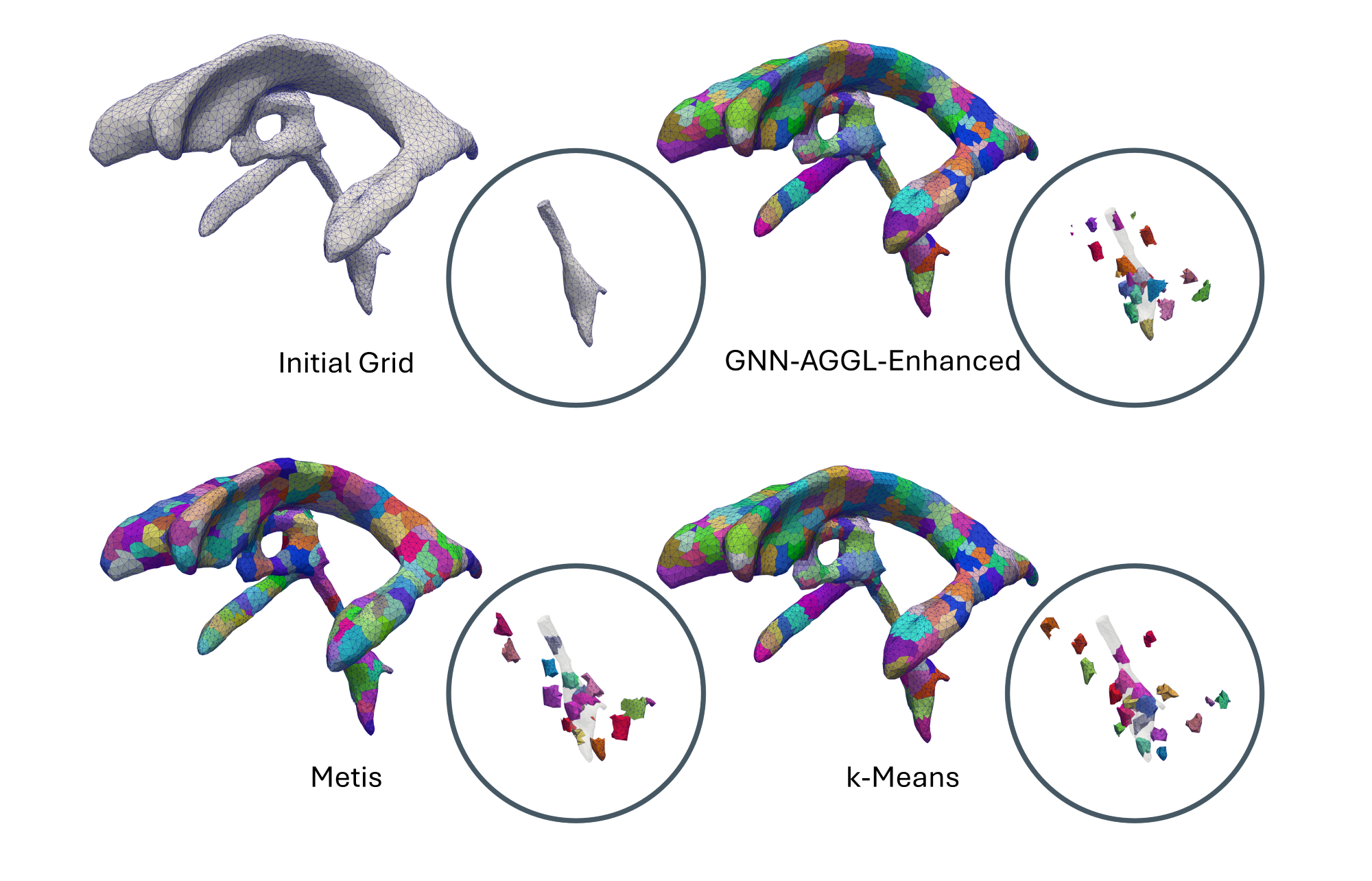}
    \caption{Test case 2: initial mesh of the human brain's ventricles and agglomerated meshes using different strategies (GNN-AGGL-Enhanced, METIS, and k-means). In the circles, we depict a detail of the exploded agglomerated grids of the fourth ventricle and aqueduct of Sylvius.}
\label{fig:agglomeratedventricles}
\end{figure}

To further test the generalization capabilities of the developed algorithm, we apply it to a more complex domain. In particular, we consider human brain ventricles reconstructed starting from a Magnetic Resonance Imaging (MRI) scan taken from the OASIS database \cite{lamontagne_oasis-3_2019}. Brain ventricles are a network of interconnected cavities filled with cerebrospinal fluid. The domain is highly non-convex and presents many convoluted pathways, constrictions, and narrowed sections. However, due to the actual research on brain fluid dynamics, constructing high-quality meshes is of major interest \cite{corti_numerical_2023,fumagalli_polytopal_2024}. The initial tetrahedral mesh consists of $28\,706$ tetrahedral elements, as shown in Figure~\ref{fig:agglomeratedventricles}. 
We agglomerated such a grid using GNN-AGGL-Enhanced, METIS, and k-means algorithms. The target mesh size of Algorithm~\ref{alg:agglalgo} used for this test case is $h/10$, where $h$ is the diameter of the initial mesh. This choice leads to agglomerated meshes with approximately 300 elements. In Figure~\ref{fig:agglomeratedventricles} we report the results of the agglomeration processes, also considering a detail of the agglomerated mesh of the fourth ventricle and aqueduct of Sylvius.  
\begin{table}[t]
\centering
\begin{tabular}{|c|Q|Q|}
\hline
\textbf{Agglomeration} & \textbf{Circle} & \textbf{Uniformity } \\
\textbf{Method} & \textbf{Ratio} & \textbf{Factor} \\ 
\hline
\textbf{GNN-AGGL-Enhanced} 
& \cellcolor[RGB]{255,114,0} 0.2252 & \cellcolor[RGB]{184,255,0} 0.6393 \\ \hline
\textbf{METIS}
& \cellcolor[RGB]{255,113,0} 0.2227 & \cellcolor[RGB]{196,255,0} 0.6151 \\ \hline
\textbf{k-Means} 
& \cellcolor[RGB]{255,118,0} 0.2313 & \cellcolor[RGB]{191,255,0} 0.6261 \\ \hline
\end{tabular}
\caption{Test case 2: computed quality metrics (CR and UF) for the agglomerated grids reported in Figure~\ref{fig:agglomeratedventricles}, obtained with different agglomeration strategies (GNN-AGGL-Enhanced, METIS, k-Means).  The colors are scaled starting from the minimum value of the metric (red) to the maximum one (green).}
\label{table:testcase2quality}
\end{table}
In Table~\ref{table:testcase2quality}, we report the quality metrics (CR and UF) computed on the agglomerated meshes. The VD metric is not reported because the starting mesh and, as a consequence, the agglomerated ones are not uniform due to the geometrical conformation of the domain. Due to this fact, the evaluation of the uniformity of the volumes of the elements is not significant for this test case. The three considered algorithms show similar results. This confirms that the proposed GNN-AGGL-Enhanced algorithm seems to be as accurate as METIS and k-means. Moreover, it indicates a good generalization capability of the GNN-AGGL-Enhanced model, considering that the computational domain is totally different from the ones included in the training set, both in terms of shape and dimensions.

\section{Numerical results: mesh agglomeration for heterogeneous domains}
\label{sec:results_het}
In this section, after introducing an additional quality metric to test the proposed agglomeration algorithm in domains featuring heterogeneous physical parameters, we introduce some test datasets to evaluate the GNN-AGGL-Enhanced model. Moreover, we present as an application, the agglomeration of a mesh representing a domain with microstructure.\\

To evaluate the regularity of the agglomerated elements, we use the quality metrics introduced in Section~\ref{sec:results_hom}. However, we need to introduce a new quantity that will be later used to measure the ability of the algorithm to preserve the different physical regions separated in the agglomeration process. We define the \textbf{Heterogeneous Elements (HE)} as the percentage of elements in the agglomerated mesh that contain discontinuities in the physical parameters.
\begin{equation}
    \operatorname{HE}=\frac{\#\mathrm{ heterogeneous\,elements}}{\# \mathrm{elements}} \%.
\end{equation}
The quantity should take values close to $0$ as much as possible, to avoid the presence of agglomerated elements containing discontinuities in the physical parameters.
\subsection{Test case 3: agglomeration capabilities assessment of test datasets}
\label{sec:test_case_3}
In this test case, we consider three different test datasets, constructed on unitary cubes, but with different distributions of physical parameters. In particular, the datasets are composed as follows:
\begin{itemize}
    \item \textbf{Dataset 1}: $25$ tetrahedral meshes with a number of elements between $373$ and $10 \ 389$ with the same physical parameter in the whole domain.
    \item \textbf{Dataset 2}: $25$ tetrahedral meshes with a number of elements between $437$ and $8 \ 549$ with two different physical parameters in two regions of the domain.
    \item \textbf{Dataset 3}: $50$ tetrahedral meshes with a number of elements between $529$ and $14 \ 960$ with two different physical parameters in more than two regions of the domain.
\end{itemize}
\begin{figure}[t]
    \centering
    \includegraphics[width=0.8\textwidth]{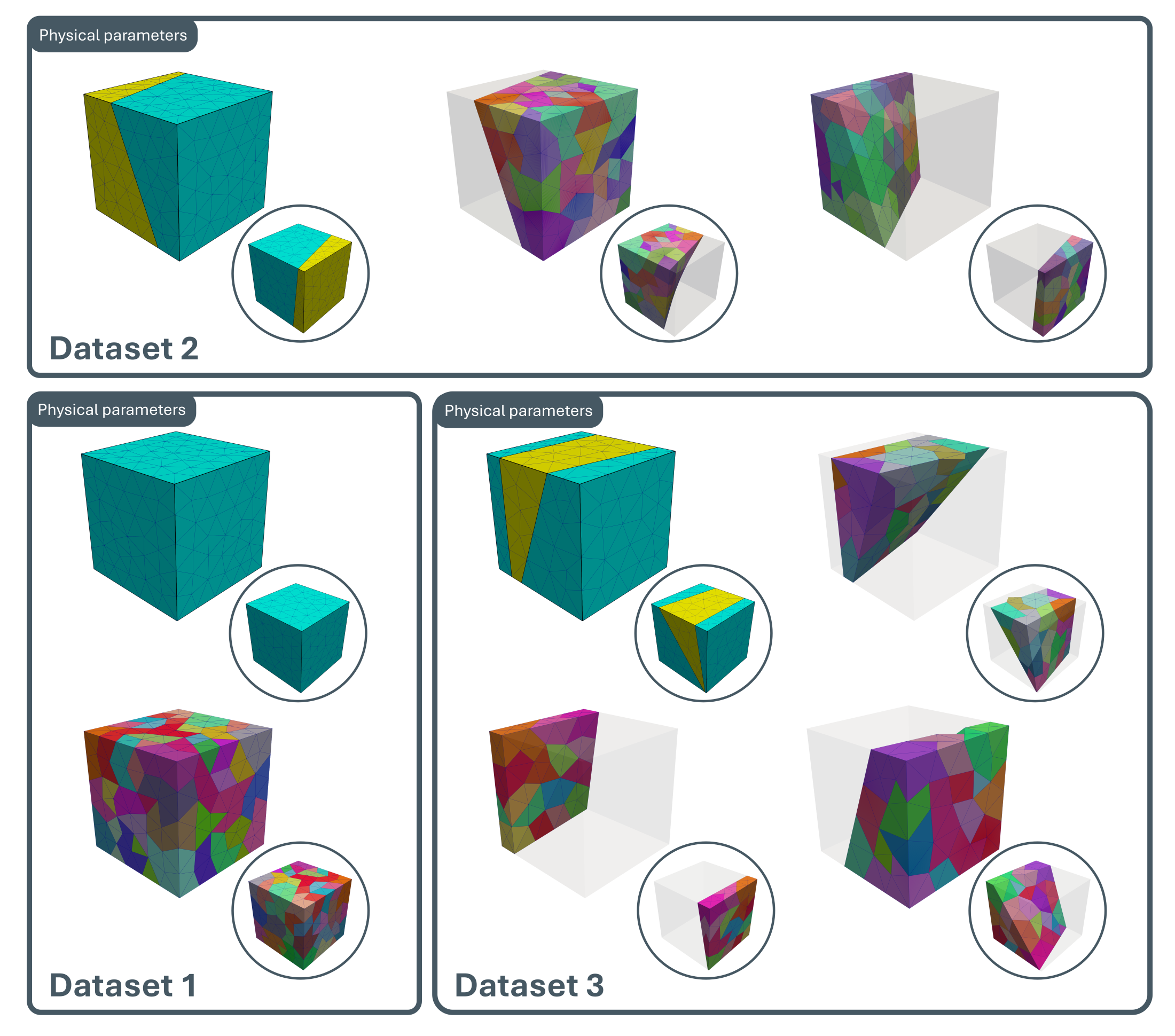}
    \caption{Test case 3: examples of agglomerated meshes, one for each dataset. The first figure (top-left) in each block represents the physical parameter heterogeneities, and the others show the agglomerated elements. Each figure is reported from two different views to visualize the agglomerated elements better.}
    \label{fig:test_case_3}
\end{figure}
The results in this section are obtained using these datasets and clustering until reaching a diameter of the agglomerated elements of $0.25$ times the diameter of the original mesh $h$. In Figure~\ref{fig:test_case_3}, we report an example of a mesh for each of the three test sets, together with the agglomerated meshes.
The first plot (top-left) in each block represents the initial mesh and the physical parameter heterogeneities, where different colors are used to represent different values of the physical parameter, and the bold lines denote the sub-region boundaries.
The other plots in each panel show the agglomerated elements. Each figure is reported from two different views to have a good visualization of the agglomerated elements.
The agglomerated meshes are visualized, separating the different physical parameter regions, to underline the level of separation of heterogeneities obtained by the GNN algorithm. In particular, only the example in Dataset 3 shows a single element containing two different physical values in the middle region. Finally, we can see that the agglomerated elements seem to be shape regular and to be relatively uniform in all the considered cases. 
\begin{table}[t]
\centering
\begin{tabular}{|q|q|q|q|q|}
\hline
\textbf{Dataset} & \textbf{CR} & \textbf{UF} & \textbf{VD} & \textbf{HE} \\ \hline
\textbf{Dataset 1} 
& \cellcolor[RGB]{255,100,0} 0.1962 
& \cellcolor[RGB]{124,255,0} 0.7568 
& \cellcolor[RGB]{173,255,0} 0.3390 
& - \\ \hline
\textbf{Dataset 2} 
& \cellcolor[RGB]{255, 98,0} 0.1918 
& \cellcolor[RGB]{123,255,0} 0.7592 
& \cellcolor[RGB]{188,255,0} 0.3692 
& \cellcolor[RGB]{  5,255,0} 1.8186\,\% \\ \hline
\textbf{Dataset 3}
& \cellcolor[RGB]{255,102,0} 0.1996 
& \cellcolor[RGB]{135,255,0} 0.7362 
& \cellcolor[RGB]{255,255,0} 0.5075 
& \cellcolor[RGB]{  7,255,0} 2.3384\,\%\\ \hline
\end{tabular}
\caption{Test case 3: average values of the computed quality metrics (CR, UF, VD, and HE) for the agglomerated grids reported in Figure~\ref{fig:test_case_3}.  The colors are scaled starting from the minimum value of the metric (red) to the maximum one (green).}
\label{table:testcase3quality}
\end{table}
To have a quantitative understanding of the agglomeration quality, we collect the average values of the quality metrics in the three different datasets in Table~\ref{table:testcase3quality}. First of all, the quantitative feedback confirms the good behavior of the proposed model in agglomerating while preserving the different physical regions. Indeed, on average, approximately only $2\%$ of the final elements are heterogeneous for the physical parameters in the starting tetrahedral elements. 
Moreover, the obtained values are coherent with what was obtained in the homogeneous case, with the values reported in Table~\ref{table:testcase1quality}. The important remark is that we do not make comparisons with either METIS or k-means here because they are not able to perform the agglomeration automatically considering the heterogeneity in the physical parameter values. 

\subsection{Test case 4: agglomeration of a mesh containing microstructures}
\begin{figure}[t]
\centering
    \begin{subfigure}[b]{0.40\textwidth}
        \resizebox{\textwidth}{!}{\includegraphics{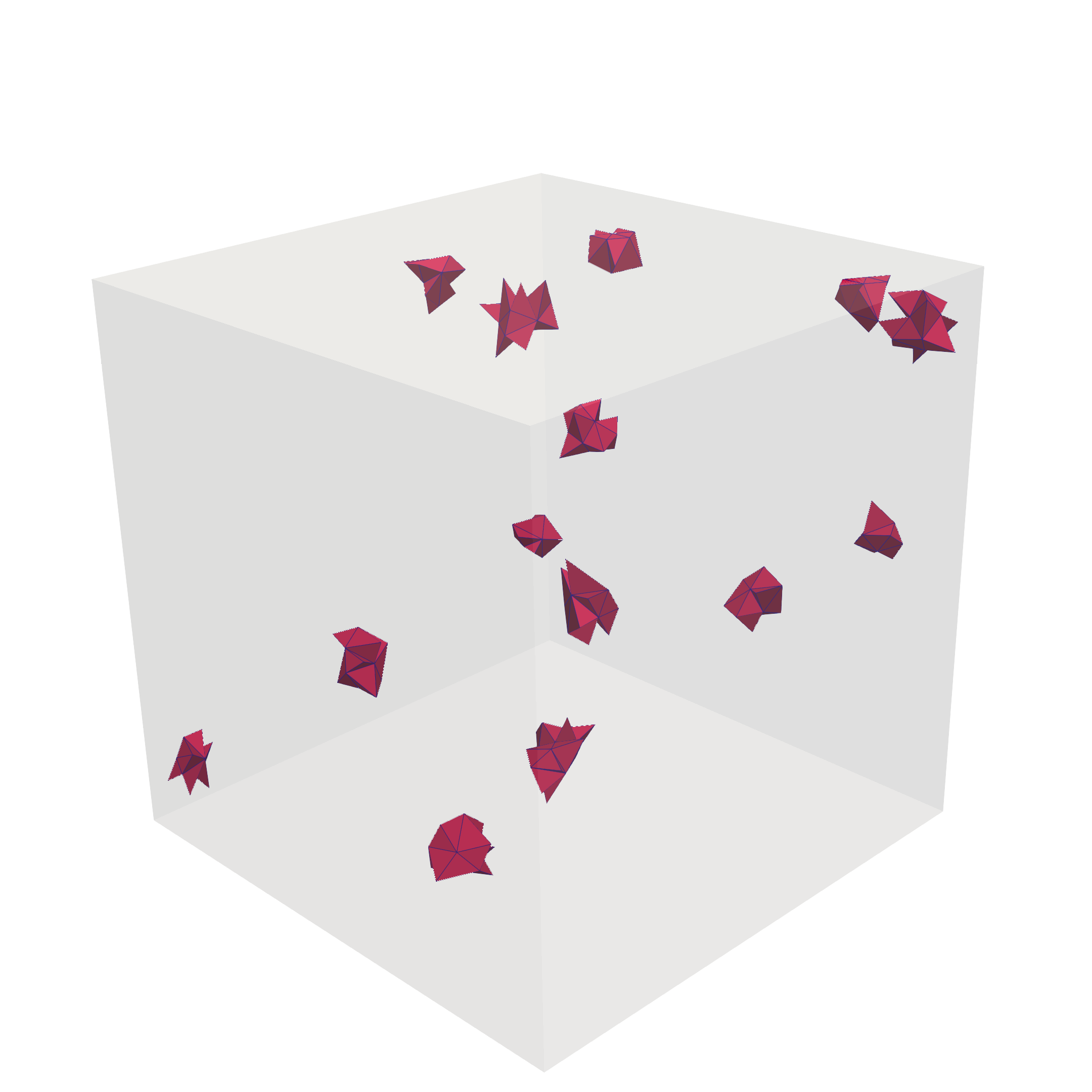}}
        \caption{Domain with microstructures}
        \label{fig:microstructures_domain}
    \end{subfigure}
    \begin{subfigure}[b]{0.40\textwidth}
    \centering
    \begin{subfigure}[b]{0.49\textwidth}
        \resizebox{\textwidth}{!}{\includegraphics{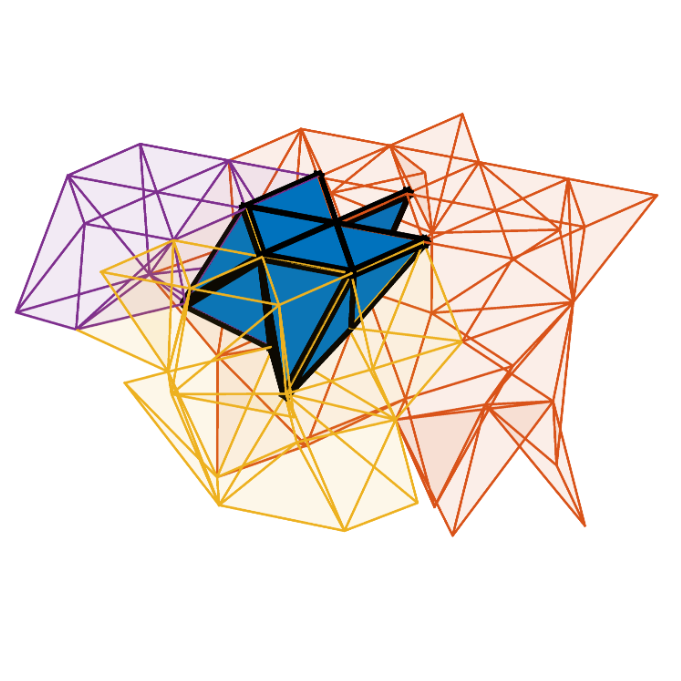}}
        \caption{Microstructure}
        \label{fig:microstructures_1}    
    \end{subfigure}     
    \begin{subfigure}[b]{0.4\textwidth}
        \resizebox{\textwidth}{!}{\includegraphics{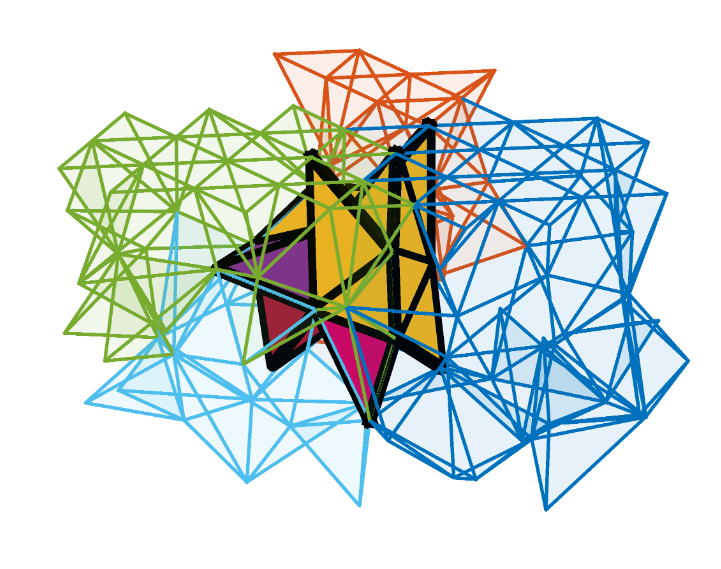}}
        \caption{Microstructure}
        \label{fig:microstructures_2}    
    \end{subfigure}      
    \begin{subfigure}[b]{0.5\textwidth}
        \resizebox{\textwidth}{!}{\includegraphics{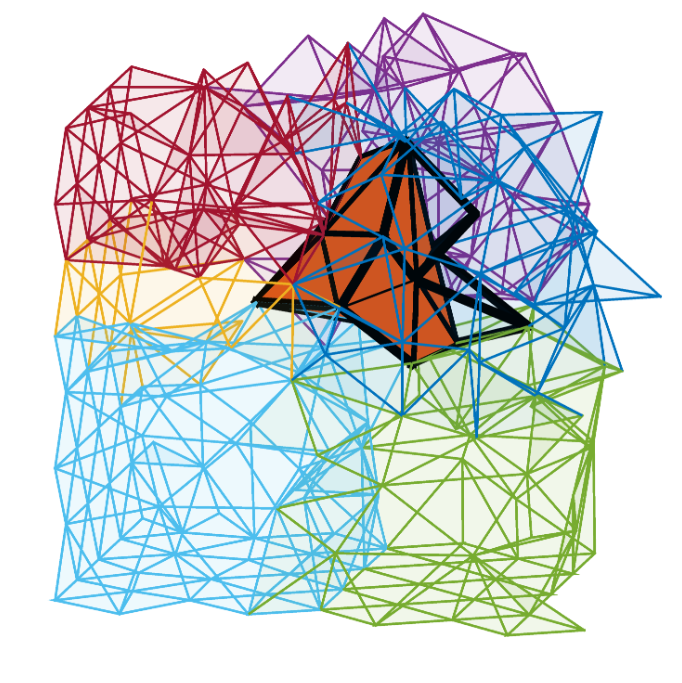}}
        \caption{Microstructure}
        \label{fig:microstructures_3}    
    \end{subfigure}  
    \end{subfigure}
    \caption{Test case 4: domain with microstructure (a), and zoom of three domains' regions containing a microstructure (b),(c), and (d).}
    \label{fig:testcase4_figures}
\end{figure}

Computational domains with complex microstructures at different scales are typical in the field of Integrated Computational Material Engineering (ICME). More specifically, many heterogeneous materials have this configuration, such as alloys, polymers, tissues, bones, and porous rocks \cite{gao_ultra_efficient_2021,trivedi_continuum_2023,maslov_biomechanical_2020,antonietti_discontinuous_2019,antonietti_polytopic_2022}.
The computational domain considered in this example is a unitary cube with $15$ microstructures in random positions, meshed by a tetrahedral grid made by $71 \ 836$ tetrahedral elements (see Figure~\ref{fig:microstructures_domain}). Each microstructure is composed of $10$ to $20$ tetrahedra, and so it is of negligible dimensions with respect to the diameter of the domain. We model this domain by setting the "physical" parameter of the background to $0$ and the microstructure one equal to $1$.
We agglomerate the underlying tetrahedral grid using the proposed algorithm until reaching a diameter of the agglomerated mesh's elements less than $0.15\,h$, where $h$ is the diameter of the original mesh. We get an agglomerated mesh of $2\,264$ elements. We observe that our model, despite the microstructures having a small size, can identify each one of them. Indeed, we found at least one agglomerated element per microstructure in the agglomerated mesh.
In Figure~\ref{fig:testcase4_figures}, each shadowed color represents a different agglomerated element, the brighter colors are the agglomerates with physical parameter $1$, and, finally, the bold black lines represent the boundary of the microstructure. We observe that the model can gather the tetrahedra of a microstructure in one agglomerated element. In some cases, the elements of a microstructure are divided into smaller agglomerates (see Figure~\ref{fig:microstructures_2}). In some instances, a few tetrahedra of a microstructure are glued in one of the adjacent background's bigger elements, as in Figure~\ref{fig:microstructures_3}. However, by computing the percentage of tetrahedra of each microstructure correctly classified in the final grid, we obtain a value around $84.90\%$.

\color{black}
\section{Robustness and generalization capabilities}
\label{sec:generalization}
In this section, we apply our algorithm to a diverse collection of meshes obtained from both heterogeneous and homogeneous domains to test its robustness and generalization capabilities. 
\subsection{Test case 5: mesh agglomeration of complex geometries}
\begin{figure}[t]
    \centering
    \includegraphics[width=\textwidth]{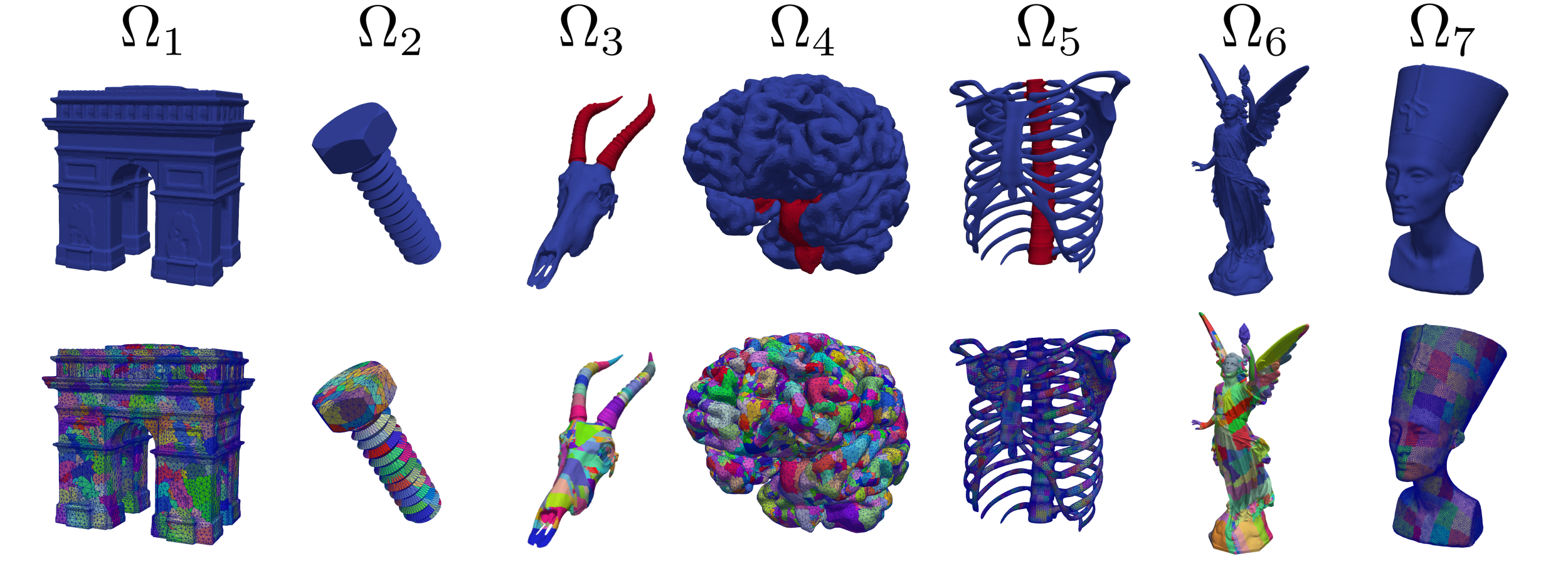}
    \caption{\color{black}Test case 5: computational meshes with indication of the heterogeneities (red and blue) for the domains $\Omega_3$,$\Omega_4$, and $\Omega_5$ (first row) and corresponding agglomerated meshes (second row). Where the mesh edges are not shown, the surface is too refined to be visualized. }
    \label{fig:test_case_5}
\end{figure}

\begin{figure}[t]
    \centering
    \includegraphics[width=\textwidth]{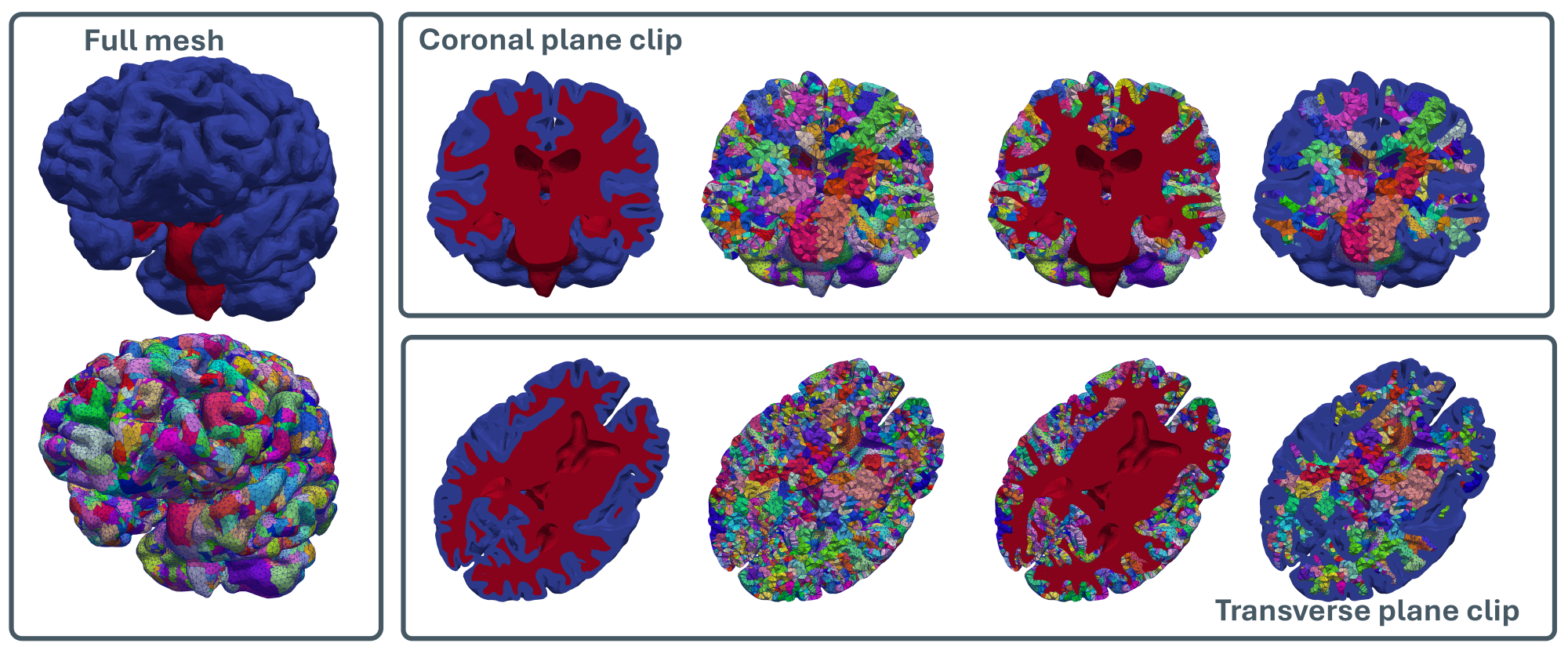}
    \caption{\color{black}Test case 5: brain geometry $\Omega_4$ with the corresponding distinction between white (in red) and grey (in blue) matters and corresponding agglomerated mesh. In the mesh visualization, the agglomerated elements are shown with different colors. We report full brain visualization, a clipped view in the coronal plane, and a transverse plane view. In the clipped views, we also show separately the agglomerated elements in the white and grey matters.}
    \label{fig:test_case_5_brain}
\end{figure}

We consider computational surfaces obtained from different open-source datasets to test the ability of our method to generalize to various mesh sizes and geometrical shapes. We report all the information on the considered examples for this test case in Table~\ref{table:testcase5meshdata}. 
Moreover, we provide a visualization of the computational domains (and corresponding fine initial tetrahedral meshes) in Figure~\ref{fig:test_case_5} (first row).
\begin{table}[t]
\color{black}
\centering
\begin{tabular}{|c|l|c|c|c|}
\hline
\multicolumn{2}{|c|}{\textbf{Domain}} & \textbf{Elements} & \textbf{Database} & \textbf{Source file type} \\ \hline
$\Omega_1$ & Arc de triomphe
& $179\,582$
& \href{https://drive.google.com/file/d/13zmGxikHiiSv9-eu8wZDTOWtPmR-KV5b/view}{FTetWild} \cite{fast_hu_2020}
& Tetrahedral mesh (model 551021)
\\ \hline
$\Omega_2$ & Bolt
& $12\,175$
& \href{https://drive.google.com/file/d/13zmGxikHiiSv9-eu8wZDTOWtPmR-KV5b/view}{FTetWild} \cite{fast_hu_2020}
& Tetrahedral mesh (model 371103)
\\ \hline
$\Omega_3$ & Damaliscus korrigum
&  $507\,480$
& \href{https://www.artec3d.com/3d-models/damaliscus-korrigum}{Artec3D}
& Surface PLY (meshed with VMTK \cite{antiga_image-based_2008}) 
\\ \hline
$\Omega_4$ & Human brain
& $609\,770$ 
& \href{https://sites.wustl.edu/oasisbrains/}{OASIS-3} \cite{lamontagne_oasis-3_2019}
& MRI Image (meshed with SVMTK \cite{mardal_mathematical_2022})
\\ \hline
$\Omega_5$ & Human rib cage
& $317\,552$
& \href{https://free3d.com/3d-model/ribcage-v2--594512.html}{Free3D}
&  Surface OBJ (meshed with VMTK \cite{antiga_image-based_2008})
\\ \hline
$\Omega_6$ & Lucy angel statue
& $179\,582$ 
& \href{https://graphics.stanford.edu/data/3Dscanrep/}{Stanford Repository}
& Surface STL (meshed with VMTK \cite{antiga_image-based_2008}) \\ \hline
$\Omega_7$ & Nefertiti bust
& $351\,086$
&  \href{https://www.thingiverse.com/thing:3974391}{Berlin Egyptian Museum}
& Surface STL (meshed with VMTK \cite{antiga_image-based_2008}) 
\\ \hline
\end{tabular}
\caption{\color{black}Test case 5: mesh information used in the test case and proper references to the databases.}
\label{table:testcase5meshdata}
\end{table}
In this test case, we consider four homogeneous domains: the Arc de Triomphe ($\Omega_1$), a bolt ($\Omega_2$), the Lucy Angel statue ($\Omega_6$), and the Nefertiti bust ($\Omega_7$), and three heterogeneous ones 
\begin{itemize} 
\item a damaliscus korrigum skull ($\Omega_3$), where we consider the difference between the bone and the horns; 
\item a human brain ($\Omega_4$) with the distinction between white and grey matter, which is fundamental in the modeling of physical phenomena \cite{corti_discontinuous_2023,corti_structure_2024};
\item a human rib cage ($\Omega_5$) with distinction of the thoracic vertebrae.
\end{itemize}
\color{black}
For all the geometries, we perform a mesh agglomeration process using the proposed GNN algorithm until reaching a diameter of the agglomerated mesh's elements less than $0.15\,\textrm{diam}(\Omega)$, where $\textrm{diam}(\Omega)$ is the diameter of the domain. The agglomerated meshes are reported in the second row of Figure~\ref{fig:test_case_5}. In the visualization, each shadowed color represents a different agglomerated element. In Table~\ref{tab:testcase5agg}, the number of elements of the agglomerated grids ($N_{\textrm{agg}}$) is compared with the corresponding number of elements of the initial tetrahedral meshes ($N_{\textrm{tets}}$). The last column of Table~\ref{tab:testcase5agg} also report the percentage reduction $\xi$ of the number of elements defined as
\begin{equation}
\label{eq:def:eta}
\xi=\dfrac{N_\mathrm{tets}-N_\mathrm{agg}}{N_\mathrm{tets}}.
\end{equation}
\begin{table}[t]
\centering
\color{black}
\begin{tabular}{|c|l|q|q|q|}
\hline
\multicolumn{2}{|c|}{\textbf{Domain}} & $N_{\textrm{agg}}$ & $N_{\textrm{tets}}$ & $\%$ reduction $\xi$ \\ \hline
$\Omega_1$ & Arc de triomphe
&  $2\,160$
&  $179\,582$
& $98,7972\%$ \\ \hline
$\Omega_2$ & Bolt
& $170$
& $12\,175$
& $98,6037\%$ \\ \hline
$\Omega_3$ & Damaliscus korrigum
& $1\,683$
& $507\,480$ 
& $99,6684\%$ \\ \hline
$\Omega_4$ & Human brain
& $22\,524$
& $609\,770$
& $96,3061\%$ \\ \hline
$\Omega_5$ & Human rib cage
& $5\,192$
& $317\,552$ 
& $98,3650\%$ \\ \hline
$\Omega_6$ & Lucy angel statue
& $1\,623$
& $179\,582$
& $99,0962\%$ \\ \hline
$\Omega_7$ & Nefertiti bust
& $900$ 
& $351\,086$
& $99,7437\%$ \\ \hline
\end{tabular}
\caption{\color{black} Test case 5: number of elements of the agglomerated grids ($N_{\textrm{agg}}$), corresponding number of elements of the initial tetrahedral grids ($N_{\textrm{tets}}$), and $\%$ reduction $\xi$ defined as in \eqref{eq:def:eta}.}
\label{tab:testcase5agg}
\end{table}
We observe that the agglomerates are coarser for homogeneous geometries with refined initial meshes and relatively simple boundaries (e.g., the Nefertiti bust). In contrast, the heterogeneous properties and complexity of the brain geometry result in a larger number of agglomerates. 
We remark that many of the meshes tested in this example are composed of more than $10$ times the elements of the most refined mesh in the training set ($\simeq 17\,000$ elements). In particular, the brain mesh is composed of more than $35$ times the elements of the finest mesh in the training set, confirming the scalability of the proposed method. This test case demonstrates the proposed algorithm’s robustness and generalization capabilities across a wide range of geometries and complexity levels.
In the heterogeneous cases, we can also observe robustness concerning the physical parameters. To have a more quantitative analysis of the results of our GNN agglomeration strategy, we computed the number of tetrahedra on the initial mesh that were incorrectly agglomerated into an element with a different physical parameter. For the three considered cases, we have that this value is $0.0108\%$ for $\Omega_3$, $1.440\%$ for $\Omega_4$, and $0.1458\%$ for $\Omega_5$. The result of the brain geometry $\Omega_4$ seems to be worse. Still, we have to underline that almost $32.64\%$ of the initial mesh elements have a face on the interface between white and grey matter, leading to a particularly challenging interface to be considered. The good quality of the distinction between white and grey matter can be observed in Figure~\ref{fig:test_case_5_brain}. These results are of particular significance, considering that neither the k-means nor METIS can automatically address this task and that the brain geometry and the internal interface are highly complex.

\color{black}
\section{Performance of the GNN agglomeration algorithm in PDE solvers}
\label{sec:solver}
\color{black}
In this section, we test the performance of the proposed agglomeration strategy when applied within PDE solvers. Specifically, we consider a three-dimensional linear elasticity equation (with unitary Lamé coefficients) with Dirichlet boundary conditions, which we discretize using a polytopal discontinuous Galerkin method.
Our model problem reads: 
\begin{equation}
\label{eq:elasticity}
\begin{dcases}
    -\nabla \cdot (2 \boldsymbol{\varepsilon}(\boldsymbol{u}) +  (\nabla\cdot\boldsymbol{u})\mathbf{I}) = \boldsymbol{f} & \mathrm{in}\,\Omega, \\
    \boldsymbol{u} = \boldsymbol{g} & \mathrm{on}\,\partial\Omega,
\end{dcases}
\end{equation}
\color{black}
where $\boldsymbol{\varepsilon}(\boldsymbol{u})$ is the symmetric part of the gradient of the displacement $\boldsymbol{u}\in \mathbf{H}^1(\Omega)$, $\boldsymbol{f}\in \mathbf{L}^2(\Omega)$ is the given forcing term, and $\boldsymbol{g}\in \mathbf{H}^{1/2}(\partial\Omega)$ is the Dirichlet boundary datum.
The discretization of problem \eqref{eq:elasticity} is based on a polytopal discontinuous Galerkin method, cf. \cite{antonietti_highorder_2018} for the method's formulation and analysis of the corresponding time-dependent problem. The discrete approximation space is defined as $\boldsymbol{V}_h^\ell = \{\boldsymbol{v} \in \mathbf{L}^2(\Omega) : \boldsymbol{v}|_K \in [\mathbb{P}^\ell(K)]^3 \, \forall K \in \mathcal{T}_h\}$, where $\ell\geq1$ is the polynomial approximation order. In the DG formulation, we choose the stabilization parameter to be equal to $10$. The subsequent numerical tests are performed considering the following manufactured  solution: 
\begin{equation*}
    \boldsymbol{u}(x,y,z) = \left[\cos(\pi x)\sin(\pi y)\sin(\pi z),\sin(\pi x)\cos(\pi y)\sin(\pi z),\sin(\pi x)\sin(\pi y)\cos(\pi z) \right]^T.
\end{equation*}
The forcing term $\boldsymbol{f}$ and the boundary condition $\boldsymbol{g}$ are computed accordingly. We compute the discretization errors as the difference between the exact solution $\boldsymbol{u}\in\mathbf{H}^1(\Omega)$ and the numerical solution $\boldsymbol{u}_h\in\boldsymbol{V}_h^\ell$ in both the $\mathbf{L}^2$-norm $\|\cdot\|_{\mathbf{L}^2}$ and the DG-norm, defined as $\|\boldsymbol{\cdot}\|_\mathrm{DG}^2 = \sum_{K\in\mathcal{T}_h} \|\boldsymbol{\varepsilon}(\boldsymbol{\cdot})\|_{\mathbf{L}^2(K)}^2 + \sum_{F\in\mathcal{F}_h} \|\eta^{1/2}\jump{\boldsymbol{\cdot}}\|_{\mathbf{L}^2(F)}^2$, where $\jump{\cdot}$ denotes the jump operator on the skeleton $\mathcal{F}_h$ of the mesh, cf. \cite{antonietti_highorder_2018} 
\subsection{Test case 6: convergence test on tetrahedral and agglomerated grids}
\begin{figure}[t]
    \begin{subfigure}[b]{0.27\textwidth}
        {\includegraphics[width=\textwidth]{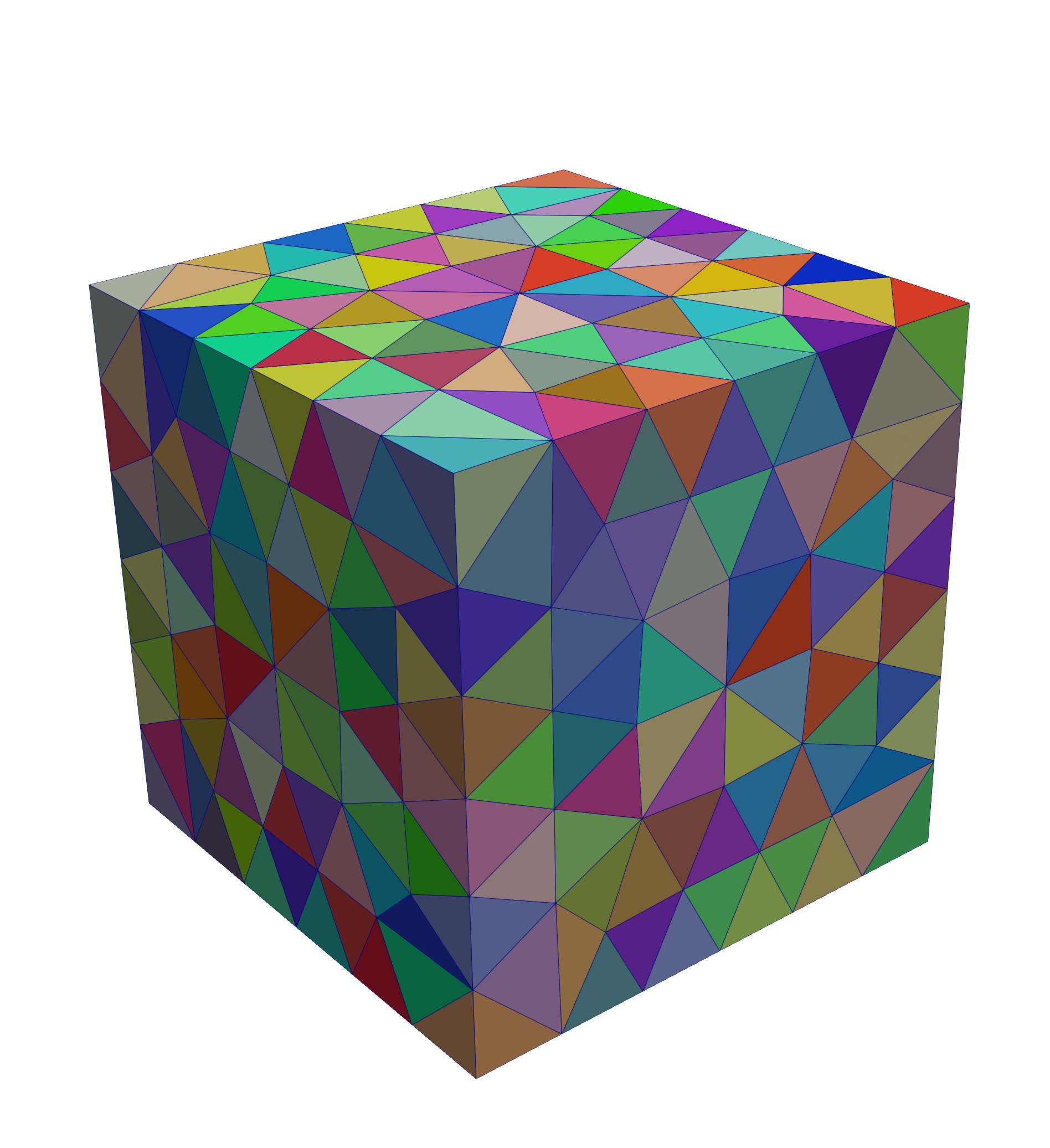}}
        \caption{\color{black}Tetrahedral mesh $\mathcal{T}_h^1$.}
        \label{fig:tc6_tetra}
    \end{subfigure}%
    \begin{subfigure}[b]{0.27\textwidth}
        {\includegraphics[width=\textwidth]{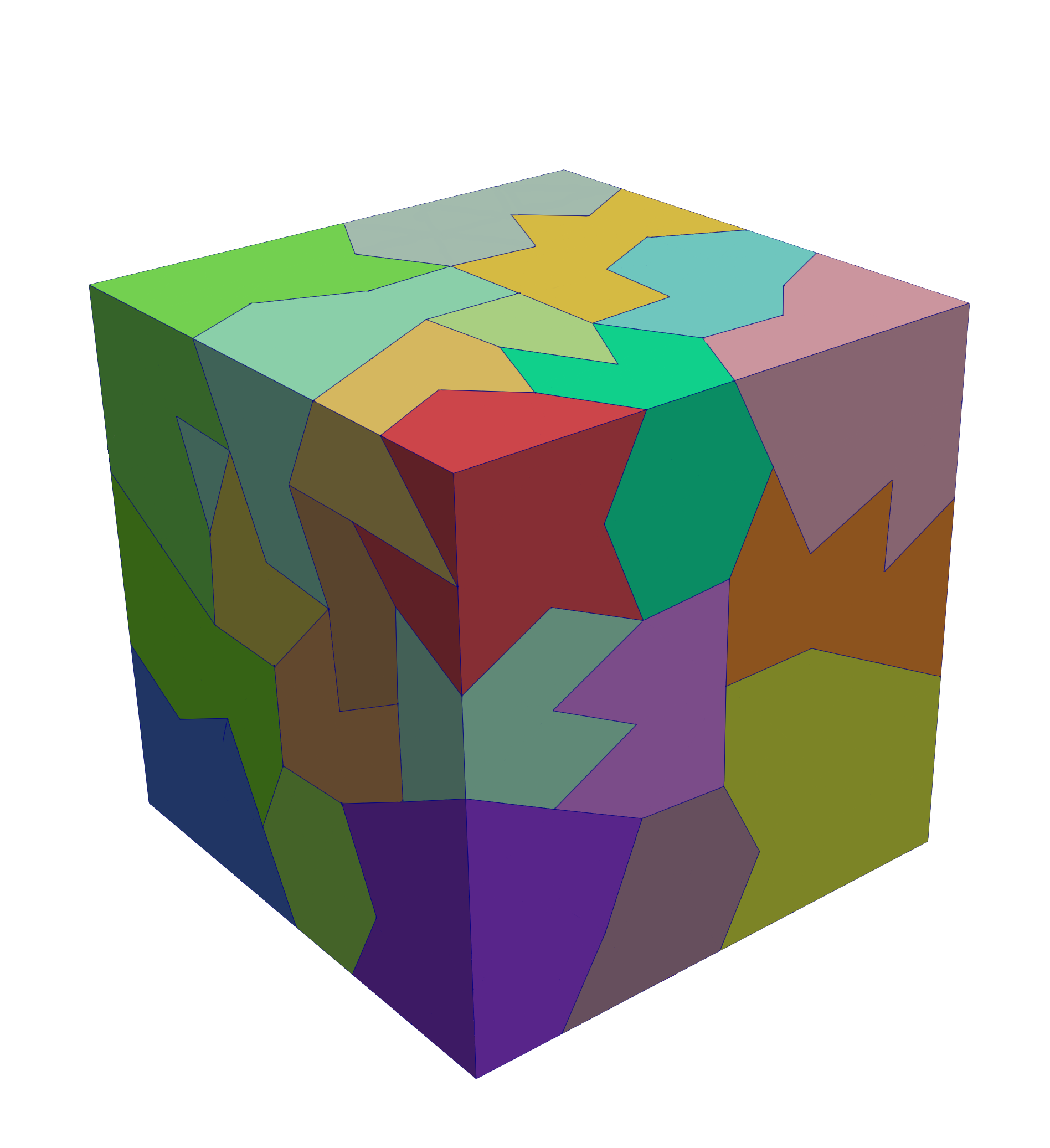}}
        \caption{\color{black}Agglomerated mesh $\mathcal{T}_h^2$.}
        \label{fig:tc6_poly}
    \end{subfigure}
    \begin{subfigure}[b]{0.45\textwidth}
          \resizebox{\textwidth}{!}{\begin{tikzpicture}

\begin{axis}[%
width=3.875in,
height=2.25in,
at={(2.6in,1.099in)},
scale only axis,
xmin=0,
xmax=60,
xminorticks=true,
xlabel = {$\mathrm{DOFs}^{1/3}$},
ylabel = {Errors},
xtick={0, 10, 20, 30, 40, 50, 60},
ymode=log,
ymin=1e-5,
ymax=10,
yminorticks=true,
axis background/.style={fill=white},
title style={font=\bfseries},
xmajorgrids,
xminorgrids,
ymajorgrids,
yminorgrids,
legend style={at={(1,1)},legend cell align=left, draw=white!15!black}
]
              
\addplot [color=blue, line width=2.0pt, mark=square*, dotted, mark options=solid]
  table[row sep=crcr]{%
25.3402  1.0797e-00 \\			
34.3919  2.3575e-01 \\
43.3318  2.1944e-02 \\
};
\addlegendentry{DG error - $\mathcal{T}_h^1$}

\addplot [color=blue, line width=2.0pt, mark=*, mark options=solid]
  table[row sep=crcr]{%
25.3402  1.5309e-02 \\			
34.3919  1.2220e-03 \\
43.3318  7.7971e-05 \\
};
\addlegendentry{$\mathbf{L}^2$ error -$\mathcal{T}_h^1$}

\addplot [color=magenta, line width=2.0pt, mark=square*, dotted, mark options=solid]
  table[row sep=crcr]{%
7.7631   2.4048e-00 \\
10.5362  1.6356e-00 \\
13.2755  6.1217e-01 \\
17.0416  1.8332e-01 \\
18.7115  3.7979e-02 \\
};

\addlegendentry{DG error - $\mathcal{T}_h^2$}

\addplot [color=magenta, line width=2.0pt, mark=*, mark options=solid]
  table[row sep=crcr]{%
7.7631   1.3591e-01 \\
10.5362  6.3531e-02 \\
13.2755  1.3110e-02 \\
17.0416  2.6317e-03 \\
18.7115  3.9114e-04 \\
};
\addlegendentry{$\mathbf{L}^2$ error -  $\mathcal{T}_h^2$}

\node[right, align=left, text=black, font=\footnotesize] at (axis cs:12,1.75) {$\ell=2$};

\node[right, align=left, text=black, font=\footnotesize] at (axis cs:6.25,5) {$\ell=1$};

\node[right, align=left, text=black, font=\footnotesize] at (axis cs:14.5,6e-1) {$\ell=3$};

\node[right, align=left, text=black, font=\footnotesize] at (axis cs:17.75,2e-1) {$\ell=4$};

\node[right, align=left, text=black, font=\footnotesize] at (axis cs:19.5,3.5e-2) {$\ell=5$};

\node[right, align=left, text=black, font=\footnotesize] at (axis cs:6.25,3e-1) {$\ell=1$};

\node[right, align=left, text=black, font=\footnotesize] at (axis cs:12,7e-2) {$\ell=2$};

\node[right, align=left, text=black, font=\footnotesize] at (axis cs:14.5,1.2e-2) {$\ell=3$};

\node[right, align=left, text=black, font=\footnotesize] at (axis cs:17.75,2.5e-3) {$\ell=4$};

\node[right, align=left, text=black, font=\footnotesize] at (axis cs:19.5,3.5e-4) {$\ell=5$};

\node[right, align=left, text=black, font=\footnotesize] at (axis cs:26.5,1.5) {$\ell=1$};

\node[right, align=left, text=black, font=\footnotesize] at (axis cs:26.5,1.5e-2) {$\ell=1$};

\node[right, align=left, text=black, font=\footnotesize] at (axis cs:35,3e-1) {$\ell=2$};

\node[right, align=left, text=black, font=\footnotesize] at (axis cs:35,1.5e-3) {$\ell=2$};

\node[right, align=left, text=black, font=\footnotesize] at (axis cs:44,2.5e-2) {$\ell=3$};

\node[right, align=left, text=black, font=\footnotesize] at (axis cs:44,7.5e-5) {$\ell=3$};

\end{axis}
\end{tikzpicture}%
}
         \caption{\color{black}Computed errors.}
        \label{fig:tc6_convergence}
    \end{subfigure}%
    \caption{\color{black}Test case 6: (a) initial tetrahedral mesh; (b) corresponding agglomerated mesh;  (c) computed errors measured in the $\mathbf{L}^2$- and DG-norms versus the cubic root of the total number of degrees of freedom  $\textrm{DOFs}$ (loglog scale).}
    \label{fig:tc6}
\end{figure}

We consider a cubic domain $\Omega = (0,1)^3$ discretized with a tetrahedral mesh $\mathcal{T}_h^1$ consisting of $1\,356$ elements (see Figure~\ref{fig:tc6_tetra}. Using our GNN agglomeration algorithm and setting as desired mesh size $0.5 \textrm{diam}(\Omega)$, we construct a polyhedral mesh $\mathcal{T}_h^2$ of $39$ elements  (see Figure~\ref{fig:tc6_poly}. We solve problem \eqref{eq:elasticity} on the mesh $\mathcal{T}_h^1$ and varying the polynomial approximation order $\ell = 1,2,3$. We repeated the same simulation on the mesh $\mathcal{T}_h^2$ with $\ell = 1,...,5$. 
The computed approximation errors for the two simulations are reported in Figure~\ref{fig:tc6_convergence}, where the errors are plotted versus the cubic root of the total number of degrees of freedom $\textrm{DoFs}$ (log-log scale). We observe that, for the considered test case—where the underlying manufactured solution is smooth—the approximation on the agglomerated mesh yields better results than that on the initial grid at a fixed level of accuracy. 
Moreover, this test case confirms that the proposed agglomeration strategy can also be effectively used in combination with a PDE solver.

\subsection{Test case 7: agglomeration of meshes with badly-shaped elements}
\begin{figure}[t]
\centering
    \begin{subfigure}[b]{0.30\textwidth}
        {\includegraphics[width=\textwidth]{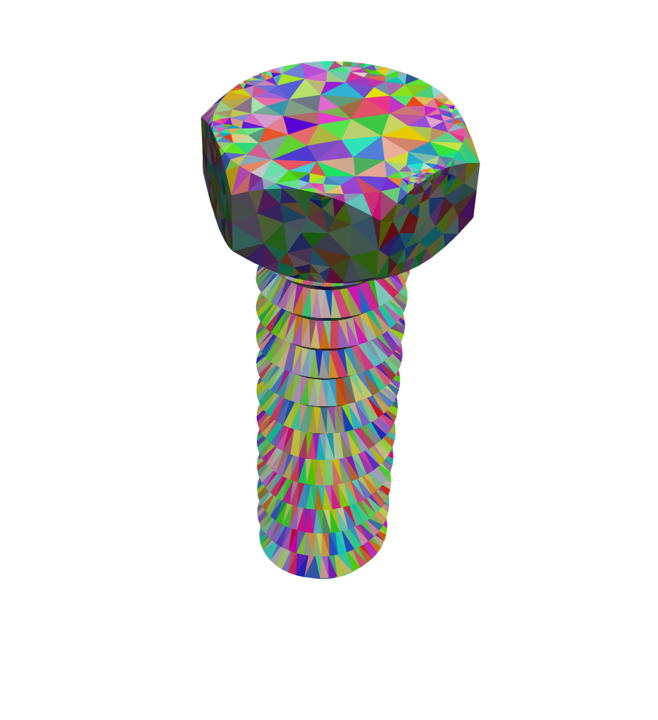}}
        \caption{\color{black}Tetrahedral mesh.}
        \label{fig:tc7_tetra}
    \end{subfigure}%
    \begin{subfigure}[b]{0.30\textwidth}
        {\includegraphics[width=\textwidth]{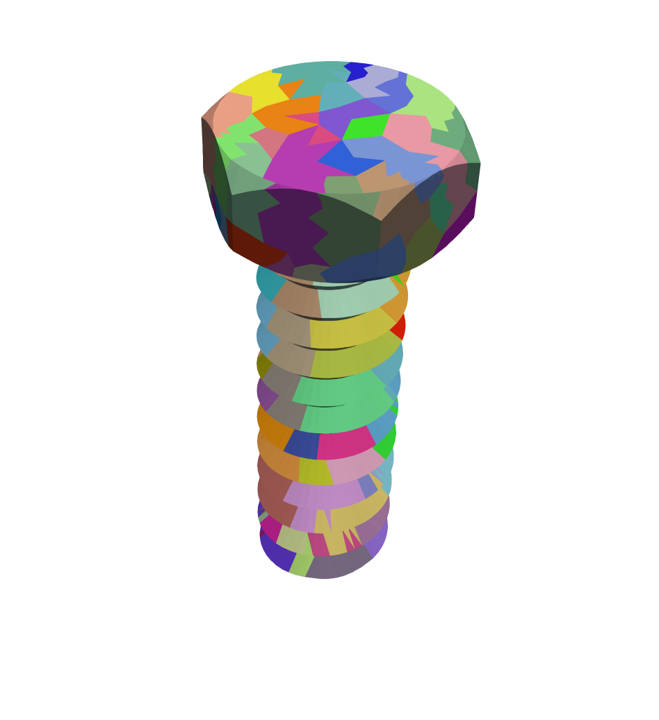}}
        \caption{\color{black}Agglomerated mesh.}
        \label{fig:tc7_poly}
    \end{subfigure}
    \begin{subfigure}[b]{0.30\textwidth}
        {\includegraphics[width=\textwidth]{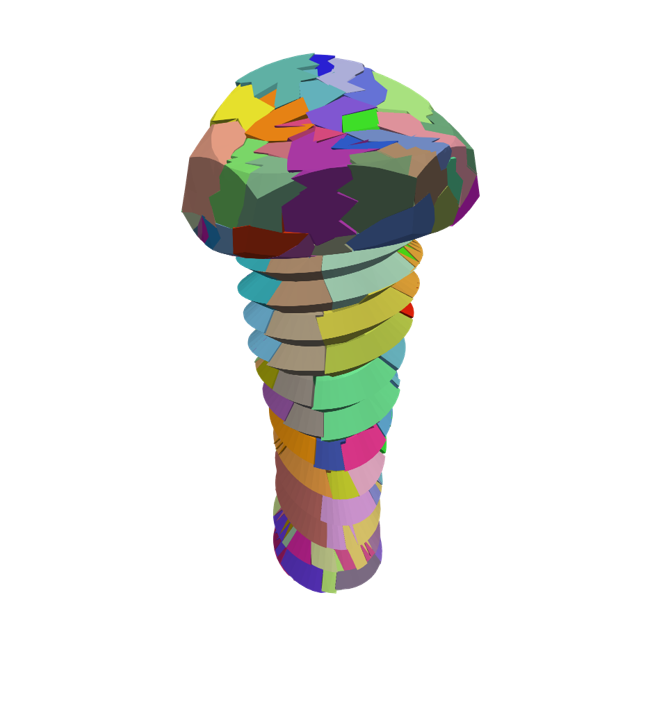}}
        \caption{\color{black} Mesh deformed by $\boldsymbol{u}_h$ for $\ell = 1$.}
        \label{fig:tc7_def}
    \end{subfigure}%
    \caption{\color{black}Test case 6: (a) initial tetrahedral mesh; (b) corresponding agglomerated mesh; (c) agglomerated mesh deformed according to the computed numerical solution for $\ell = 1$.}
    \label{fig:tc7}
\end{figure}

The last test case considers the solution of the linear elasticity problem \eqref{eq:elasticity} on the bolt domain $\Omega_2$ considered in Section~\ref{sec:generalization}. 
The initial tetrahedral mesh, consisting of 12,175 elements, was downloaded from the FTetWild database \cite{fast_hu_2020} (see Figure~\ref{fig:tc7_poly}). We note that this mesh contains severely distorted elements, and the goal of this test is to demonstrate that the proposed GNN agglomeration algorithm can also improve overall mesh quality. Using our algorithm, we construct a polyhedral mesh with 170 elements (see Figure~\ref{fig:tc7_poly}). We solve problem \eqref{eq:elasticity} on both the initial tetrahedral mesh with $\ell = 1$ and the agglomerated mesh for $\ell = 1, \dots, 4$. Figure~\ref{fig:tc7_def} reports the deformed mesh according to the computed numerical solution for $\ell = 1$.
The computed results align well with what is expected for the agglomerated mesh, confirming the effectiveness of the proposed GNN algorithm when used within a PDE solver. We also compute the approximation errors in both the $\mathbf{L}^2$ and DG norms; the results are reported in Table~\ref{tab:testcase7}.
On the agglomerated mesh, we obtain accurate results for all tested values of $\ell$. In contrast, the initial tetrahedral mesh fails to produce a stable numerical solution due to the presence of severely distorted elements, which negatively affect the solver’s accuracy. These results suggest that the agglomeration process not only enables accurate solutions but also improves mesh quality by reducing the occurrence of highly stretched tetrahedra.
\begin{table}[t]
\centering
\begin{tabular}{|c|c|q|q|}
\multicolumn{4}{c}{\textbf{Computed errors on the agglomerated mesh}} \\ \hline
\textbf{Polynomial Degree} & \textbf{Number of DOFs} & $\mathbf{L}^2-$\textbf{Error} & \textbf{DG}$-$\textbf{Error} \\ \hline
$\ell=1$ & $2\,040$
& $1.54\times 10^{-2}$
& $7.94\times 10^{-1}$ \\ \hline
$\ell=2$ & $5\,100$
& $4.08\times 10^{-3}$
& $2.69\times 10^{-1}$ \\ \hline
$\ell=3$ & $10\,200$
& $5.44\times 10^{-4}$
& $5.74\times 10^{-2}$ \\ \hline
$\ell=4$ & $17\,850$
& $5.85\times 10^{-5}$
& $8.59\times 10^{-3}$ \\ \hline
\multicolumn{4}{c}{\textbf{Computed errors on the tetrahedral mesh}} \\ \hline
\textbf{Polynomial Degree} & \textbf{Number of DOFs} & $\mathbf{L}^2-$\textbf{Error} & \textbf{DG}$-$\textbf{Error} \\ \hline
$\ell=1$ & $146\,1000$
& $3.44\times 10^{-1}$
& $1.61\times 10^{+2}$ \\ \hline
\end{tabular}
\caption{\textcolor{black}{Test case 7: computed approximation errors for the agglomerated mesh with $\ell=1,...,4$ and the tetrahedral mesh with $\ell=1$.}}
\label{tab:testcase7}
\color{black}
\end{table}

\color{black}
\section{Conclusions}
\label{sec:conclusions}
\color{black}
In this work, we presented a novel agglomeration algorithm for three-dimensional grids using a graph partitioning approach based on Graph Neural Networks (GNNs). Our model, which extends the results of \cite{antonietti_agglomeration_2024}, learns the geometric structure of mesh elements automatically and efficiently, while remaining independent of specific geometries or underlying differential models. This flexibility enables the method to handle complex scenarios—such as heterogeneous domains and microstructures—that classical strategies cannot process automatically, while preserving the geometric quality of the original mesh and significantly reducing computational costs. 
We tested the proposed approach across diverse test cases, including real geometries reconstructed from medical images and domains with microstructures and/or featuring complex geometric details. We compared the results with the analogous ones obtained by employing graph partitioners such as METIS and k-means on various test cases. Results show that our method maintains high mesh quality while efficiently handling both the geometry, the mesh topology, and the underlying physical parameters. Additionally, the GNN model offers lower online computational cost. The proposed approach has also been tested when used within PDE solvers, confirming its robustness and generalization capabilities.

Future work will focus on exploring reinforcement learning frameworks, as developed by \cite{gatti_graph_2022}, as well as extending the benchmarking considering domain-specific partitioning strategies such as SCOTCH \cite{pellegrini:hal-00770422,chevalier:hal-00402893} or Coarsening Schemes for Graph Partitioning (e.g. \cite{Ron2011407, Safro2012369})
\color{black}

\section*{Acknowledgments}
The brain MRI images were provided by OASIS-3: Longitudinal Multimodal Neuroimaging: Principal Investigators: T. Benzinger, D. Marcus, J. Morris; NIH P30 AG066444, P50 AG00561, P30 NS09857781, P01 AG026276, P01 AG003991, R01 AG043434, UL1 TR000448, R01 EB009352. AV-45 doses were provided by Avid Radiopharmaceuticals, a wholly-owned subsidiary of Eli Lilly. \textcolor{black}{We acknowledge Ilario Mazzieri (MOX - Dipartimento di Matematica, Politecnico di Milano, Italy) for the helpful discussions.}

\section*{Declaration of competing interests}
The authors declare that they have no known competing financial interests or personal relationships that could have appeared to influence the work reported in this article.
\bibliographystyle{hieeetr}
\bibliography{cas-refs}

\end{document}